\documentclass[11pt,a4paper]{article}
\usepackage[T1]{fontenc}
\usepackage[utf8]{inputenc}
\usepackage{lmodern}
\usepackage[margin=2.5cm]{geometry}
\usepackage{amsmath,amssymb,amsthm}
\usepackage{graphicx}
\usepackage{booktabs}
\usepackage{xcolor}
\usepackage[numbers,square,sort&compress]{natbib}
\usepackage[colorlinks=true,linkcolor=blue!50!black,citecolor=blue!50!black,urlcolor=blue!50!black]{hyperref}
\hypersetup{pdftitle={Round-Based Contagion Generates Chaos from Oscillatory Dynamics: A Two-Stage Contagion Map},pdfauthor={Francisco J. Mu\~noz, Juan Carlos Nu\~no, Luca Meacci}}
\usepackage{caption}
\usepackage{authblk}

\newtheorem{proposition}{Proposition}
\newtheorem{definition}{Definition}
\theoremstyle{remark}
\newtheorem{remark}{Remark}

\newcommand{\Rt}{\mathcal{R}_\tau}
\newcommand{\rT}{\rho_{T}}
\newcommand{\rNS}{\rho_{NS}}
\newcommand{\pT}{p_{T}}
\newcommand{\pNS}{p_{NS}}

\title{Round-Based Contagion Generates Chaos from Oscillatory Dynamics: A Two-Stage Contagion Map}
\author[1]{Francisco J. Mu\~noz\thanks{Corresponding author: \texttt{fcomunoz@alumnos.upm.es}}}
\author[1]{Juan Carlos Nu\~no}
\author[2,3]{Luca Meacci}
\affil[1]{\small Departamento de Matem\'atica Aplicada, Universidad Polit\'ecnica de Madrid, 28040 Madrid, Spain}
\affil[2]{\small Dipartimento di Matematica e Informatica ``Ulisse Dini'', Universit\`a degli Studi di Firenze, Viale Morgagni, 67/a, 50134, Firenze, Italy}
\affil[3]{\small Dipartimento di Ingegneria dell'Informazione (DINFO), Universit\`a degli Studi di Firenze, Via di S.~Marta, 3, 50139, Firenze, Italy}
\date{}

\begin{document}
\maketitle

\begin{abstract}
Many social processes unfold in rounds, such as elections, meetings or turns of a game, in which individuals update synchronously. We show that this round structure can turn a regular oscillation into chaos. In Katriel's differential equations, two-stage contagion (adoption after two contacts with adopters) combined with population turnover produces bistability and relaxation oscillations, but no chaos. From the elementary reactions we derive a deterministic map $F_\tau$ for a large population updated in rounds of length $\tau$: each individual makes at most one transition per round, events compete as exponential risks, and adoptions take effect from the next round. Katriel's equations are recovered as $\tau\to0$. For short rounds the map follows the relaxation cycle, and the phase of successive outbreaks obeys a circle homeomorphism, which excludes chaos. As the rounds lengthen, this circle map loses invertibility and the invariant curve breaks down; once a round lasts about two mean adoption times, comparable to the duration of an outbreak, the attractor is chaotic, while Katriel's equations remain on a periodic relaxation cycle. With the probability $p$ of adoption after a single contact as bifurcation parameter, the map goes from extinction through a window of chaos to a stable equilibrium; for $\beta/\gamma=46$, $\mu/\gamma=0.002$, $\gamma\tau=5.5$ the window is $0.023\le p\le0.056$, bounded by the invasion threshold and by a boundary crisis, whereas Katriel's equations oscillate for $0.022<p<0.41$. Chaos is supported by the Lyapunov spectrum (Kaplan--Yorke dimension $1.81$), the correlation dimension, the 0-1 test, ordinal-pattern complexity, saddle periodic orbits and a horseshoe of the phase map. One-stage contagion, $p=0$ and permanent adoption showed no chaos.
\end{abstract}

\noindent{\raggedright\textbf{Keywords:} Round-based contagion; two-stage contagion; relaxation oscillations; torus breakdown; circle map; chaos.\par}

%======================================================================
\section{Introduction}

Many social processes are organised in rounds. Votes are cast at elections, positions are revised at periodic meetings, players change strategy between the turns of a game, and choices made during one round become visible to others only at the next. Contagion in such processes is round-based: individuals update synchronously, once per round, using the state of the population at the end of the previous round. Here we ask what this round structure does to a contagion process whose continuous-time description oscillates, and show that it can turn a regular oscillation into chaos.

Many social contagions require reinforcement: an individual adopts a behaviour, opinion or technology only after repeated exposure~\cite{Granovetter1978,Watts2002,CentolaMacy2007,DoddsWatts2004}. Katriel~\cite{Katriel2019} formulated the simplest compartmental version of this idea. A naive susceptible ($S_0$) becomes primed ($S_1$) after a first contact with an adopter ($A$) and adopts after a second contact; a fraction $p$ of first contacts leads to adoption directly. Adopters eventually stop adopting ($R$), and the population is renewed by turnover. The parameter $p$ interpolates between pure two-stage contagion ($p=0$) and ordinary one-stage contagion ($p=1$). The resulting system of ordinary differential equations (ODE) shows backward bifurcations, bistability, critical transitions, Hopf bifurcations, relaxation oscillations and excitability. Heidecke and Barbarossa~\cite{HeideckeBarbarossa2021} found comparable complexity in a related model, and spatial versions show that bistability and oscillations depend on the range of interaction~\cite{Meacci2025}. Chaos, however, has not been reported in these continuous-time models. Their relaxation oscillations make them a natural test case for the question above.

For round-based processes a discrete-time description is not an approximation but the natural model. Discrete-time contagion models, however, are usually written down directly or obtained by discretising a differential equation~\cite{Allen1994,Hu2014}, so that the time step has no meaning of its own. Here, instead, we derive the model from first principles. We write the two-stage mechanism as a reaction scheme, state explicitly how individuals update in rounds, and obtain a discrete-time stochastic model whose transition probabilities follow from these rules alone. The construction is a binomial version of the $\tau$-leap method of stochastic chemical kinetics~\cite{Gillespie2001,ChatterjeeVlachos2005}, with one modification that keeps the expected exposure to adopters consistent for all step sizes. Katriel's ODE is recovered in the joint limit of small round length and large population, so that the length $\tau$ of a round is the only parameter separating the two descriptions. Throughout this paper we study the deterministic large-population limit of the model, a three-dimensional map $F_\tau$; the effects of demographic noise in finite populations are left for future work.

Our central finding can be stated in one sentence: \emph{a continuous-time relaxation oscillator becomes chaotic when updates are forced to occur in synchronous rounds whose duration is comparable to the outbreak duration.} Katriel's ODE supplies the oscillator: a slow recharge of naive individuals by turnover, lasting hundreds of mean adoption times, followed by a fast outbreak that lasts about one mean adoption time $1/\gamma$. As long as the rounds are short compared with an outbreak, the map follows the cycle and its attractor is a smooth invariant curve. When a round lasts about two adoption times, each outbreak is compressed into one or two rounds, its size depends sensitively on where the threshold crossing falls within the sequence of rounds, and the cycle breaks up into a chaotic attractor. The ratio of the round length to the outbreak duration, and not the period of the cycle, is what separates regular from chaotic dynamics.

Our main results are the following.
\begin{enumerate}
\item Round-based updating turns oscillation into chaos. For fixed rates of Katriel's model, the map follows the relaxation cycle when the rounds are short (an invariant closed curve carrying quasi-periodic or mode-locked motion) and becomes chaotic when a round lasts about two mean adoption times, through the break-up of the invariant curve, which we identify with the loss of invertibility of the circle map that links the phases of successive outbreaks; chaos ends when the turnover per round $\mu\tau$ exceeds about $0.04$. Throughout the chaotic window Katriel's ODE is on a relaxation cycle: the rounds turn an oscillation that already exists in the continuous model into chaos.
\item The two-stage parameter $p$ controls where this happens: with all rates and the round length fixed, chaos occupies a window of small $p$ bounded below by the invasion threshold $\pT$ and above by a boundary crisis just beyond a subcritical Neimark--Sacker point $\pNS$, whereas the ODE oscillates over a range more than ten times wider. One-stage contagion ($p=1$) is never chaotic.
\item Within the window, in all our tests, the chaotic attractor is unique up to coexistence with the stable equilibrium or with periodic orbits; it has a single positive Lyapunov exponent and an outbreak--recharge geometry that does not change with $p$.
\item Chaos is confirmed by independent indicators: besides a positive Lyapunov exponent it is detected by the 0-1 test, a broadband spectrum, ordinal-pattern complexity and a correlation dimension that matches the Kaplan--Yorke dimension; it is robust to precision, coordinates, initial conditions and run length; its periodic orbits and unstable manifold are consistent with a horseshoe; and along the round length the phase map of successive outbreaks has a horseshoe whose entropy bound is close to the Lyapunov exponent.
\item We are not aware of previous reports of chaos in two-stage or multi-stage contagion models, and the route to it (break-up of an invariant curve followed by a crisis) differs from the period-doubling route usual in discrete-time epidemic maps.
\end{enumerate}

Section~\ref{sec:related} reviews related discrete and multi-stage contagion models. Section~\ref{sec:elem} gives the reaction scheme; Section~\ref{sec:model} states the modelling principles, interprets the round length, defines the discrete model and derives its limits; Section~\ref{sec:local} analyses equilibria and local bifurcations; and Section~\ref{sec:chaos} shows how the rounds turn the oscillation into chaos: how chaos emerges as the rounds lengthen, the role of $p$ as bifurcation parameter, the time evolution and independent indicators of chaos, what Katriel's model does in the chaotic window, the mechanism and conditions for chaos, and the robustness of the results.

%======================================================================
\section{Related models}\label{sec:related}

Among discrete-time models, the closest antecedent to ours is the SWIR model~\cite{Choi2018}, in which a susceptible in contact with an infected individual becomes either infected or weakened ($W$), and a weakened individual in contact with an infected one becomes infected, each with its own probability. These are exactly our reactions $R_1$--$R_3$ (Section~\ref{sec:elem}), with $W$ playing the role of $S_1$. The SWIR model is defined in discrete time steps and is analysed through a nonlinear map, but it has no renewal of susceptibles: a single outbreak runs to an absorbing state, and the questions addressed are the hybrid percolation transition and the size dependence of the latent period before the explosive outbreak on random networks. The generalised epidemic process of Janssen, M\"uller and Stenull~\cite{Janssen2004} introduces the same weakened intermediate state, through the reactions $S+X\to Y+X$ and $Y+X\to2X$, in a spatial setting, and produces tricritical and discontinuous percolation transitions~\cite{Bizhani2012}; cooperative co-infection~\cite{Chen2013} shares the reinforcement mechanism. These are likewise single-outbreak processes without renewal of susceptibles.

A second family of models represents reinforcement through memory, thresholds or stages. In the generalised contagion model of Dodds and Watts~\cite{DoddsWatts2004,DoddsWatts2005}, individuals contact one another in discrete time steps, accumulate doses over a memory window and adopt when a threshold is exceeded; recovered individuals may become susceptible again. The type of collective behaviour depends only on the memory length and on the probabilities of adoption after one and after two exposures, a role analogous to that of $p$ here. Threshold models~\cite{Granovetter1978,Watts2002} and their multi-stage generalisation~\cite{Melnik2013}, in which adopters at different stages exert different influence, show that one additional stage already produces cascades impossible in single-stage models. Higher-order (simplicial) contagion~\cite{Iacopini2019}, also studied in discrete time through microscopic Markov chains~\cite{Gomez2010,Chowdhary2021}, produces explosive transitions and bistability. In the settings studied, these models reach absorbing or stationary states; to our knowledge, chaos has not been reported in any of them.

Chaos has, by contrast, been reported in discrete-time one-stage epidemic models. Allen~\cite{Allen1994} showed that discrete SIS models, and SI and SIR models with births, can show period doubling and chaos, whereas without feedback to the susceptible class they behave like their continuous analogues as long as solutions stay positive. Later work analysed the flip bifurcation and computed Lyapunov exponents in such maps~\cite{Hu2014}. In these models chaos arises through period doubling of a logistic-type map when the step is large. Many such maps are obtained by forward-Euler-type discretisation of an ODE. It has also been shown that discrete-time and continuous-time contagion dynamics coincide as the step tends to zero but can differ substantially for finite steps~\cite{Fennell2016}.

In continuous time, Katriel~\cite{Katriel2019} and Heidecke and Barbarossa~\cite{HeideckeBarbarossa2021} found that two-stage contagion combined with a source of renewal of susceptibles (demographic turnover, or waning and fading) gives rise to bistability, oscillations and complex bifurcations, and conjectured that this coupling is the origin of the rich dynamics. Our model contains Katriel's as its continuous limit, and the chaos we report occurs where his model oscillates.

The model studied here differs from all of the above in combining the two-stage mechanism, renewal of susceptibles and synchronous updating, and in being derived from the elementary processes rather than postulated or obtained from a numerical scheme. It is chaotic, its chaos depends on the two-stage parameter $p$, and the route to chaos is the break-up of an invariant curve followed by a boundary crisis, not period doubling.

%======================================================================
\section{Elementary processes}\label{sec:elem}

Consider a well-mixed population of $N$ individuals in the classes $S_0$, $S_1$, $A$, $R$. The two-stage mechanism with turnover consists of the following elementary processes, written as reactions with rate constants per unit time:
\begin{align*}
&R_1:\ S_0+A\xrightarrow{\ p\beta/N\ }2A, &&R_2:\ S_0+A\xrightarrow{\ (1-p)\beta/N\ }S_1+A, &&R_3:\ S_1+A\xrightarrow{\ \beta/N\ }2A,\\
&R_4:\ A\xrightarrow{\ \gamma\ }R, &&R_{5,6,7}:\ X\xrightarrow{\ \mu\ }S_0\quad(X=S_1,A,R).
\end{align*}
$R_1$--$R_3$ are contagion events catalysed by adopters: direct adoption, priming and adoption of primed individuals. $R_4$ is cessation of adoption, and $R_5$--$R_7$ describe turnover at constant population: every individual is replaced by a naive one at rate $\mu$. In chemical terms, $R_1$--$R_7$ describe a two-step autocatalysis in a continuous-flow stirred reactor. Mass-action kinetics for the fractions $s_0=S_0/N$, $s_1$, $a$, $r$ gives Katriel's temporary-adoption model~\cite{Katriel2019}:
\begin{equation}\label{eq:ode}
\dot s_0=\mu-\beta s_0a-\mu s_0,\quad \dot s_1=(1-p)\beta s_0a-\beta s_1a-\mu s_1,\quad \dot a=\beta(ps_0+s_1)a-(\gamma+\mu)a,
\end{equation}
with $r=1-s_0-s_1-a$. Throughout, time is measured in units of the mean adoption time $1/\gamma$, and we use the dimensionless parameters
\[
\rho=\beta/\gamma,\qquad \nu=\mu/\gamma,\qquad p,
\]
together with the dimensionless round length $\gamma\tau$, which we write simply as $\tau$ when no confusion arises. In Eq.~\eqref{eq:ode}, adoption invades the adoption-free state $E_0=(1,0,0)$ if and only if $\rho p/(1+\nu)>1$. Katriel~\cite{Katriel2019} uses instead the contact parameter $\delta=\beta/(\gamma+\mu)=\rho/(1+\nu)$ and the duration ratio $\gamma/\mu=1/\nu$, so that his invasion threshold $\delta=1/p$ is the same condition. Our reference point $\rho=46$, $\nu=0.002$ corresponds to $\delta\approx45.9$ and $\gamma/\mu=500$: adoption periods much shorter than the residence time in the population, strong contact and small $p$. This regime lies outside the one illustrated in his figures (for instance $\delta\approx4$--$8$, $\gamma/\mu=10$--$20$, $p=0.15$), so all comparisons with his model below use Eq.~\eqref{eq:ode} with our rates.

%======================================================================
\section{The discrete-time model}\label{sec:model}

\subsection{Modelling principles and the meaning of the round}\label{sec:principles}

The discrete model rests on five principles. The first two are those of the reaction scheme itself, the next two specify how time is organised, and the last one is a closure assumption on the exposure within a round.
\begin{enumerate}
\item[(P1)] \emph{Mass action.} Each adopter makes contacts with randomly chosen members of the population at rate $\beta$, and each contact with a naive or primed individual triggers $R_1$, $R_2$ or $R_3$ with the probabilities of the reaction scheme.
\item[(P2)] \emph{Memoryless elementary events.} Cessation, turnover and contagion are independent exponential risks; when several can occur, the first one determines the transition (competing risks).
\item[(P3)] \emph{Rounds.} Time is divided into rounds of length $\tau$. In each round an individual can change class at most once, and a change of class made during a round becomes effective, and visible to others, from the next round.
\item[(P4)] \emph{Conditional independence.} Given the state at the start of a round, individuals evolve independently during the round.
\item[(P5)] \emph{Exposure closure.} Within a round, the total contagion hazard generated by the adopters is spread uniformly over the round.
\end{enumerate}
Principles (P1), (P2) and (P4) are the standard assumptions behind the chemical master equation, and hence behind Katriel's ODE. Principle (P3) is the essential additional ingredient, and it is what distinguishes processes that unfold in rounds from those that unfold in continuous time. Principle (P5) is a closure assumption that yields closed-form transition probabilities; (P1)--(P4) alone determine the exact exposure profile, and Remark~\ref{rem:profile} shows that the main results are robust to this choice.

Many social changes happen in rounds of this kind. In elections, voters are exposed to campaigns and to the opinions of others between consecutive elections, but their choice is registered only at the next election, and each voter makes at most one choice per election. In iterated games, and in tournaments or periodic negotiations, players observe the strategies used in one turn and revise their own for the next. Committees, assemblies and periodic meetings revise positions at each meeting. Fashions and consumer choices follow collections, seasons or release cycles. School and workplace practices change by terms or by course. On online platforms, behaviour is often organised in daily or weekly cycles of exposure and response. In all these cases the time between rounds is set by the institution or the activity, not by the speed of contagion, and an individual who changes behaviour during a round affects others mainly from the next round on.

Within the model, the round is characterised by two dimensionless numbers. The product $\gamma\tau$ compares the round with the mean time $1/\gamma$ during which an adopter actively promotes the behaviour; numerically, we find chaos only when $\gamma\tau\gtrsim2$ (about $1.7$ at the lowest in our scans), i.e.\ when this active phase is shorter than half a round, as for short-lived enthusiasm, campaigns, trends or bursts of mobilisation. The product $\mu\tau$ is the fraction of the population renewed per round; in our scans chaos also requires $\mu\tau\lesssim0.04$, i.e.\ members who stay on average for at least some 25 rounds, as in a weekly or monthly cycle within a community whose members remain for years. Processes with long rounds, such as elections every few years, share the round structure but would fall outside this range unless membership is very stable; the chaotic regime is therefore most relevant for processes with frequent rounds. The parameter $p$ measures how often a single exposure suffices for adoption, and $1-p$ how often an intermediate stage of awareness or priming is needed.

\subsection{Definition}

Fix a round length $\tau>0$ and write $X_n=(S_{0,n},S_{1,n},A_n,R_n)\in\mathbb{N}^4$, with $S_{0,n}+S_{1,n}+A_n+R_n=N$, for the counts at time $n\tau$. Principles (P1)--(P5) translate into four rules for the round from $n\tau$ to $(n+1)\tau$.

\begin{enumerate}
\item[(i)] \emph{At most one transition per individual and round} (P2, P3). Individuals change class at most once in each round, according to the first of their competing elementary events.
\item[(ii)] \emph{Adopters act from the next round} (P3). The $A_n$ adopters present at the beginning of the round are the only sources of contagion during the round. Adopter $i$ is active for a time $T_i=\min(E_i,\tau)$, where $E_i\sim\mathrm{Exp}(\gamma+\mu)$ is the time of its first event ($R_4$ or $R_6$). Individuals that become adopters during the round start acting at the next one.
\item[(iii)] \emph{Contagion} (P1, P2, P5). Given the active times, each susceptible is exposed to a contagion hazard that is constant over the round (P5) and has integral $\Lambda_n=(\beta/N)\sum_{i\le A_n}T_i$, the total adopter activity during the round. Turnover competes with contagion as an independent exponential risk of rate $\mu$.
\item[(iv)] \emph{Independence} (P4). Given $X_n$ and the active times, individuals evolve independently.
\end{enumerate}

These rules turn the reaction scheme into explicit transition probabilities. For a susceptible, contagion and turnover are competing risks with integrated hazards $\Lambda_n$ and $\mu\tau$, so that the probability of some event in the round is $P_n=1-e^{-(\Lambda_n+\mu\tau)}$ and the event is contagion with probability $\phi_n=\Lambda_n/(\Lambda_n+\mu\tau)$. Writing $\kappa=1-e^{-(\gamma+\mu)\tau}$ for the probability that an adopter has an event within the round, the transitions are:

\begin{center}\small
\begin{tabular}{@{}lll@{}}
\toprule
class at $n\tau$ & transitions in $[n\tau,(n+1)\tau)$ & probability\\
\midrule
$S_0$ & $\to A$ ($R_1$), $\to S_1$ ($R_2$) & $p\,\phi_nP_n$, \ $(1-p)\,\phi_nP_n$\\
$S_1$ & $\to A$ ($R_3$), $\to S_0$ ($R_5$) & $\phi_nP_n$, \ $(1-\phi_n)P_n$\\
$A$ & $\to R$ ($R_4$), $\to S_0$ ($R_6$) & $\tfrac{\gamma}{\gamma+\mu}\kappa$, \ $\tfrac{\mu}{\gamma+\mu}\kappa$\\
$R$ & $\to S_0$ ($R_7$) & $1-e^{-\mu\tau}$\\
\bottomrule
\end{tabular}
\end{center}

\noindent Turnover of a naive individual leaves it naive. The counts at $(n+1)\tau$ are obtained from multinomial draws for each class with these probabilities, so that $(X_n)$ is a Markov chain that conserves the population exactly. No step of this construction refers to Eq.~\eqref{eq:ode}: the probabilities follow from (P1)--(P5). We use the chain only as the microscopic foundation of the deterministic map studied below.

\begin{remark}[time profile of the exposure]\label{rem:profile}
Principles (P1)--(P4) determine the exposure profile itself: given the active times, the contagion hazard at time $t$ within the round is $(\beta/N)\,\#\{i:T_i>t\}$, which decreases as adopters stop. Replacing it by a constant hazard with the same integral (P5) leaves unchanged the expected number of contacts and the probability $1-e^{-(\Lambda_n+\mu\tau)}$ that a susceptible has some event in the round, because this probability depends only on the integrated hazard. Only the split of these events between contagion and turnover changes, by a relative amount of order $\mu\tau$. To check the robustness of our results to (P5), we repeated the main computations with the exact profile, for which the probability of contagion in the large-population map is $q_c=\int_0^\tau\lambda(t)\exp\!\big(-\int_0^t\lambda-\mu t\big)\,dt$ with $\lambda(t)=\beta a\,e^{-(\gamma+\mu)t}$, instead of $\phi P$. The largest Lyapunov exponent changed by less than $5\%$ ($0.0371\to0.0372$ per step at the reference point of Section~\ref{sec:timeevol}; $0.0186\to0.0189$ at $p=0.035$; $0.0172\to0.0181$ at $\rho=30$, $\gamma\tau=3$, $p=0.06$; $0.0080\to0.0083$ at a weakly chaotic point), the invasion threshold by less than $0.5\%$, and the chaotic window in $p$ from $[0.0230,0.0555]$ to $[0.0235,0.0565]$ at the reference rates and from $[0.0350,0.0885]$ to $[0.0350,0.0815]$ at $\rho=30$, $\gamma\tau=3$, with 89 and 90 chaotic values of $p$ respectively. We keep (P5) because it gives closed-form transition probabilities.
\end{remark}

\subsection{Limits}

\begin{proposition}[short rounds]\label{prop:small}
As $\tau\to0$, every transition probability equals the rate of the corresponding reaction in $R_1$--$R_7$, evaluated at $X_n$, times $\tau$, plus $O(\tau^2)$. Hence the chain converges to the continuous-time Markov chain with mass-action propensities, i.e.\ to the chemical master equation of the reaction scheme; for small $\tau$ it is a first-order (binomial) $\tau$-leap approximation of it~\cite{Gillespie2001,ChatterjeeVlachos2005}.
\end{proposition}

\begin{proof}
$\sum_iT_i=A_n\tau+O(\tau^2)$, so $\Lambda_n=\beta A_n\tau/N+O(\tau^2)$, $\phi_nP_n=\Lambda_n+O(\tau^2)$, $(1-\phi_n)P_n=\mu\tau+O(\tau^2)$ and $\kappa=(\gamma+\mu)\tau+O(\tau^2)$.
\end{proof}

For $N\to\infty$ at fixed $\tau$, $\sum_iT_i/N\to a\,\mathbb{E}[T]=a\,(1-e^{-(1+\nu)\tau})/(\gamma(1+\nu))$ by the law of large numbers, with $a=A_n/N$. The fractions $x_n=X_n/N$ then follow, over any finite number of rounds, the deterministic map defined next. This is an exact large-population limit of the chain, not a truncation: no expansion in $\tau$ is involved.

\begin{definition}
For $x=(s_0,s_1,a)$ and $r=1-s_0-s_1-a$, let
\[
\Lambda=\frac{\rho\,(1-e^{-(1+\nu)\tau})}{1+\nu}\,a,\quad m=\nu\tau,\quad P=1-e^{-(\Lambda+m)},\quad \phi=\frac{\Lambda}{\Lambda+m},\quad \kappa=1-e^{-(1+\nu)\tau},
\]
with $\tau$ measured in units of $1/\gamma$. The map $F_\tau:(s_0,s_1,a)\mapsto(s_0',s_1',a')$ is
\begin{align}
s_0'&=s_0(1-\phi P)+s_1(1-\phi)P+a\,\tfrac{\nu}{1+\nu}\,\kappa+r\,(1-e^{-m}),\nonumber\\
s_1'&=(1-p)\,s_0\,\phi P+s_1(1-P),\label{eq:map}\\
a'&=p\,s_0\,\phi P+s_1\,\phi P+a\,(1-\kappa).\nonumber
\end{align}
\end{definition}

\begin{proposition}[consistency]\label{prop:consistency}
$F_\tau$ maps the simplex $\Delta=\{s_0,s_1,a\ge0,\ s_0+s_1+a\le1\}$ into itself, and $F_\tau(x)=x+\tau f(x)+O(\tau^2)$ uniformly on $\Delta$, where $f$ is the vector field~\eqref{eq:ode} in units of $1/\gamma$. Thus, for small $\tau$, $F_\tau$ behaves as a smooth one-step method of order one for Katriel's equations.
\end{proposition}

\begin{proof}
Invariance holds because the right-hand sides are expectations of a probability vector. Expanding in $\tau$: $\Lambda=\rho a\tau+O(\tau^2)$, $\phi P=\Lambda+O(\tau^2)$, $(1-\phi)P=\nu\tau+O(\tau^2)$, $\kappa=(1+\nu)\tau+O(\tau^2)$, which gives $s_0'=s_0+\tau(\nu-\rho s_0a-\nu s_0)+O(\tau^2)$ and the analogous expansions of $s_1'$ and $a'$.
\end{proof}

Proposition~\ref{prop:consistency} is a statement about the limit $\tau\to0$. For finite $\tau$, $F_\tau$ is not a numerical scheme for Eq.~\eqref{eq:ode} but the mean-field dynamics of a process organised in rounds: it preserves positivity and the population for every $\tau$, its transition probabilities never exceed one, and its parameters keep their individual-level meaning. A forward Euler step of Eq.~\eqref{eq:ode} with $\gamma\tau=3$, by contrast, can produce negative fractions once $\rho a\tau>1$.

\begin{remark}[why integrated exposure]
The standard $\tau$-leap freezes the propensities at the beginning of the step, which amounts to $\Lambda=\rho\tau a$: each adopter would then act during the whole round even if, for $\gamma\tau\gtrsim1$, it most likely stops adopting early in it. With that rule, contagion per adopter grows without bound with $\tau$ and the invasion threshold drifts away from Katriel's value (for $\gamma\tau\to\infty$ it tends to $\rho p\tau\approx1$ instead of $\rho p/(1+\nu)=1$). We found chaos with the frozen rule too, but in a region where the ODE with the same rates has no endemic adoption at all, which makes the comparison with Katriel's model meaningless. Rule~(ii) removes this artefact; its only cost is that adopters created within a round act from the next one, which is the content of principle (P3).
\end{remark}

%======================================================================
\section{Equilibria and local bifurcations}\label{sec:local}

\begin{proposition}[invasion threshold]\label{prop:threshold}
The adoption-free state $E_0=(1,0,0)$ is a fixed point of $F_\tau$. Its eigenvalues are $e^{-m}$ (twice) and
\[
r_0=1+\big(1-e^{-(1+\nu)\tau}\big)\big(\Rt-1\big),\qquad \Rt=\frac{\rho p}{1+\nu}\cdot\frac{1-e^{-\nu\tau}}{\nu\tau}.
\]
Hence $E_0$ is stable if and only if $\Rt<1$. In particular, for $p=0$ it is stable for every $\rho$, and $\Rt\to\rho p/(1+\nu)$, Katriel's threshold, as $\tau\to0$.
\end{proposition}

\begin{proof}
At $E_0$, $\partial(\phi P)/\partial a=\rho(1-e^{-(1+\nu)\tau})(1-e^{-m})/((1+\nu)m)$, so $\partial a'/\partial a=p\,\partial(\phi P)/\partial a+e^{-(1+\nu)\tau}$; the other eigenvalues come from the $s_0$ and $s_1$ equations at $a=0$.
\end{proof}

The threshold $\Rt=1$ can be read in two ways. For fixed $p$ it defines $\rT(\tau)$, which differs from Katriel's value $(1+\nu)/p$ only by the factor $\nu\tau/(1-e^{-\nu\tau})\approx1+\nu\tau/2$; for $\nu=0.005$ and $p=0.05$ it increases from $20.10$ to $20.35$ as $\gamma\tau$ goes from $0$ to $5$. For fixed rates and round length it defines the critical single-contact probability
\begin{equation}\label{eq:pT}
\pT=\frac{(1+\nu)\,\nu\tau}{\rho\,(1-e^{-\nu\tau})},
\end{equation}
below which adoption cannot invade.

\begin{proposition}[endemic equilibria]\label{prop:branch}
For fixed $\Lambda$ the map~\eqref{eq:map} is affine, $x'=M(\Lambda)x+c$ with $c=(1-e^{-m},0,0)$. The endemic equilibria are therefore exactly the points
\[
x^*(\Lambda)=\big(I-M(\Lambda)\big)^{-1}c,\qquad \rho(\Lambda)=\frac{(1+\nu)\,\Lambda}{\big(1-e^{-(1+\nu)\tau}\big)\,a^*(\Lambda)},\qquad \Lambda>0,
\]
whenever $x^*(\Lambda)\in\Delta$.
\end{proposition}

This parametrisation makes continuation of the endemic branch, of its eigenvalues and of its Neimark--Sacker (NS) points straightforward. The criticality of an NS point is given by the first Lyapunov coefficient $d$ of the map~\cite{Kuznetsov}, which we evaluated with finite-difference second and third differentials; the procedure reproduces to six digits the value $d=3.2553$ obtained with exact differentials for a related map. For $p=0.05$ and $\nu=0.005$ the branch has a single NS point $\rNS(\tau)$, which is \emph{subcritical} ($d>0$) for all $\tau$ examined:

\begin{center}\small
\begin{tabular}{@{}lccccccccc@{}}
\toprule
$\gamma\tau$ & 0.05 & 0.5 & 1 & 2 & 3 & 4 & 5 & 7 & 10\\
\midrule
$\rNS$ & 90.6 & 77.8 & 63.7 & 43.9 & 32.6 & 25.8 & 21.3 & 16.0 & 11.9\\
$d$ & 0.20 & 3.97 & 7.64 & 7.51 & 5.54 & 4.12 & 3.21 & 2.18 & 1.47\\
\bottomrule
\end{tabular}
\end{center}

\noindent As $\tau\to0$, $\rNS$ tends to the Hopf point of Eq.~\eqref{eq:ode}, which for these parameters is at $\rho=91.8$, and $d\to0$, as expected for a map that approximates a time-$\tau$ flow. The endemic equilibrium is unstable for $\rT<\rho<\rNS(\tau)$; in this band the attractor is a large-amplitude oscillation or, as shown next, chaotic. The NS curve descends rapidly with $\tau$, so that the oscillatory band shrinks and closes near $\gamma\tau\approx5.3$, where $\rNS=\rT$. Being subcritical, the NS bifurcation also produces a band of bistability above $\rNS(\tau)$, where large-amplitude dynamics coexists with the stable equilibrium. In the same way, for fixed $(\rho,\nu,\tau)$ the equilibrium loses stability as $p$ decreases through an NS point $\pNS$, which we locate by bisection on the spectral radius; at the three points used in Section~\ref{sec:p} the first Lyapunov coefficient is $d=10.6$, $6.7$ and $4.5$, so these bifurcations are subcritical as well.

%======================================================================
\section{Rounds turn the oscillation into chaos}\label{sec:chaos}

We computed the largest Lyapunov exponent $\lambda_1$ of $F_\tau$ with a renormalised twin trajectory, starting from two initial conditions and using $3\times10^4$--$6\times10^4$ transient and averaging steps (protocol in Appendix~\ref{app:num}). Exponents are given per step, i.e.\ per round, unless stated otherwise. An attractor is classified as chaotic if $\lambda_1>0.003$ per step, as an equilibrium if $\log_{10}a$ varies by less than $10^{-3}$, as an oscillation otherwise, and as extinction if $a<10^{-30}$ throughout the averaging window. Because a cycle spans tens to hundreds of rounds, this threshold is conservative near the onset of chaos; where needed we also use the exponent per cycle (Section~\ref{sec:roundlength}).

\subsection{From oscillation to chaos as the rounds lengthen}\label{sec:roundlength}

\begin{figure}[t]
\centering
\includegraphics[width=\linewidth]{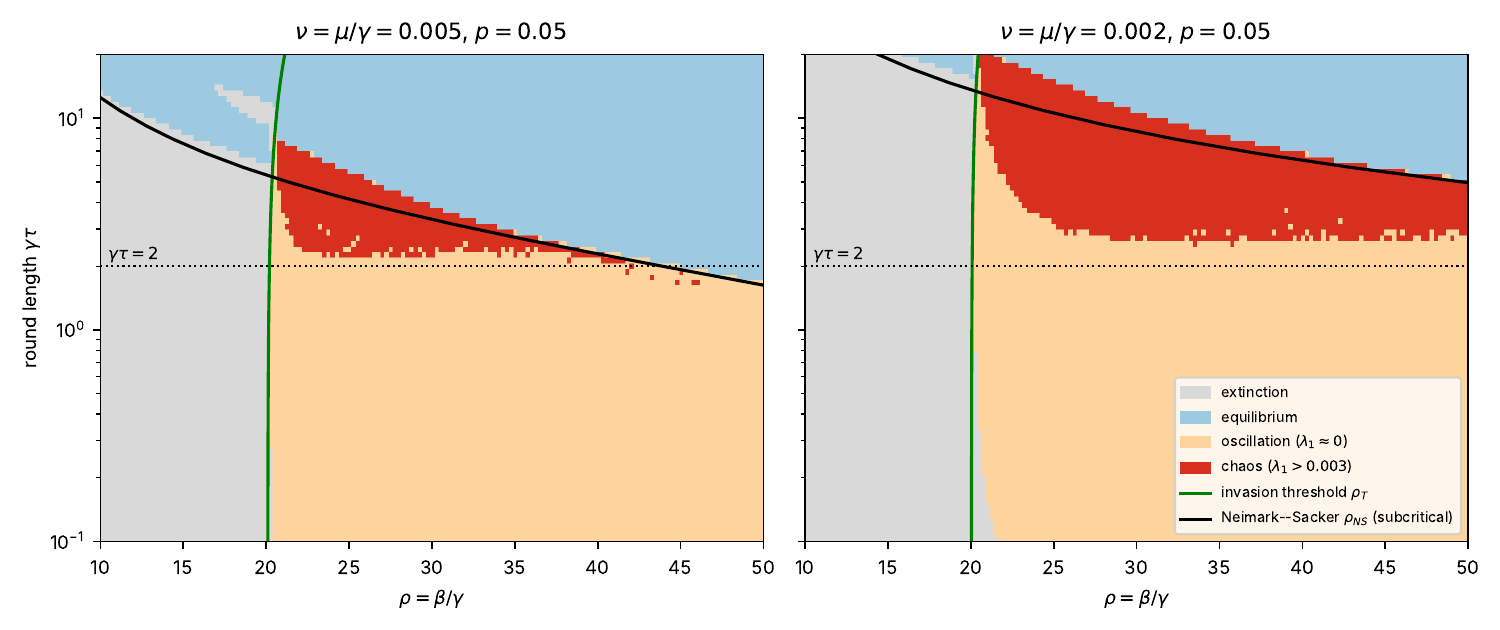}
\caption{Attractors of $F_\tau$ in the plane of the contact ratio $\rho=\beta/\gamma$ and the round length $\gamma\tau$ ($p=0.05$). The bottom of each panel is close to Katriel's continuous model. Green: invasion threshold $\rT(\tau)$ (Proposition~\ref{prop:threshold}); black: subcritical Neimark--Sacker curve $\rNS(\tau)$; dotted: $\gamma\tau=2$. Near the invasion threshold the period of the relaxation cycle can exceed the observation window, and some oscillatory points at small $\gamma\tau$ are then labelled as extinction.}
\label{fig:classes}
\end{figure}

We first follow Katriel's oscillation as the rounds lengthen. Figure~\ref{fig:classes} shows the attractors in the plane $(\rho,\gamma\tau)$ for $p=0.05$ and two turnover ratios. For short rounds, every $\rho$ above the invasion threshold gives an oscillation with $\lambda_1\approx0$: the attractor is an invariant closed curve, the discrete counterpart of Katriel's relaxation cycle (Section~\ref{sec:why}). Chaos appears when the round reaches about two mean adoption times. For $\nu=0.005$ it occupies $1.7\lesssim\gamma\tau\lesssim7.6$ and $20.8\le\rho\le46$; for $\nu=0.002$, $2.6\lesssim\gamma\tau\le20$ (the edge of the explored window) and $20.5\le\rho\le50$ (likewise the edge). Within these ranges it fills most of the band between $\rT$ and $\rNS(\tau)$, together with the bistable strip just above $\rNS(\tau)$; for $\nu=0.005$ the band closes at $\gamma\tau\approx5.3$, and the chaotic points at larger $\gamma\tau$ lie in the region where the equilibrium is stable and coexists with the chaotic attractor. The upper limit in $\tau$ corresponds in both cases to a turnover per round $\mu\tau=\nu\gamma\tau\approx0.04$: beyond it, a large fraction of the population is renewed in each round, the oscillatory band closes and the equilibrium becomes stable.

The same transition is seen along lines of fixed rates, i.e.\ for a fixed parameter set of Katriel's model with rounds of increasing length: up to $\gamma\tau\approx2$ the exponent fluctuates around zero, then becomes positive, reaching $0.005$--$0.009$ per unit time, and finally negative when the equilibrium becomes stable. Along the round length at fixed $\rho=25$, $\nu=0.005$, $p=0.05$ the route can be followed through a one-dimensional reduction (Fig.~\ref{fig:circlemap}; details in Supplementary Material, Section~\ref{Ssec:circle}). Let $\varphi_k\in[0,1)$ be the phase, measured in rounds, at which the $k$-th outbreak crosses the level $a=10^{-3}$. Transverse perturbations contract by a factor of $e^{-5}$--$e^{-9}$ per cycle (second Lyapunov exponent between $-0.011$ and $-0.021$ per step), so that successive phases obey a degree-one circle map $\varphi_{k+1}=g(\varphi_k)$, up to a scatter that grows from $3\times10^{-6}$ near $\gamma\tau=1.8$ to $4\times10^{-2}$ at $\gamma\tau=3$ (Fig.~\ref{fig:circlemap}a). For $\gamma\tau\le1.39$, $g$ is an orientation-preserving homeomorphism: its minimum slope decreases from $0.90$ at $\gamma\tau=1$ to $0.002$ at $\gamma\tau=1.39$ (Fig.~\ref{fig:circlemap}b). By the theory of circle homeomorphisms~\cite{deMeloVanStrien1993}, its rotation number is then unique, its orbits are quasi-periodic or mode-locked and its topological entropy is zero, so the attractor cannot be chaotic; this is the regime of the invariant curve. At $\gamma\tau\approx1.40$ the minimum slope of $g$ vanishes and $g$ ceases to be invertible: the invariant curve can no longer carry a circle homeomorphism, which is the breakdown of the invariant torus~\cite{Aronson1982,AfraimovichShilnikov1983}. Invertibility is recovered for $1.54\lesssim\gamma\tau\lesssim1.74$ and lost for good above $\gamma\tau=1.75$. Positive exponents occur only where $g$ is not invertible: at isolated values from $\gamma\tau=1.80$ and persistently from $\gamma\tau\approx1.98$, interspersed with periodic windows, as in circle maps beyond criticality~\cite{Jensen1984} (Fig.~\ref{fig:circlemap}c). For $\gamma\tau\gtrsim2.3$ (and at isolated values from $2.16$), $g$ contains $k\ge3$ disjoint arcs, each of which is mapped across the whole circle with a margin of at least five times the scatter; this one-dimensional horseshoe forces a topological entropy of at least $\ln k$ per cycle. At $\gamma\tau=3$, $k=11$, which gives at least $0.0076$ per step, against $\lambda_1=0.0080$. Exponents per cycle are larger than exponents per step by the number of steps per outbreak, $300$--$950$ along this line; with the per-step threshold $\lambda_1>0.003$ used in the scans, the first chaotic value appears only at $\gamma\tau=2.18$.

\begin{figure}[t]
\centering
\includegraphics[width=\linewidth]{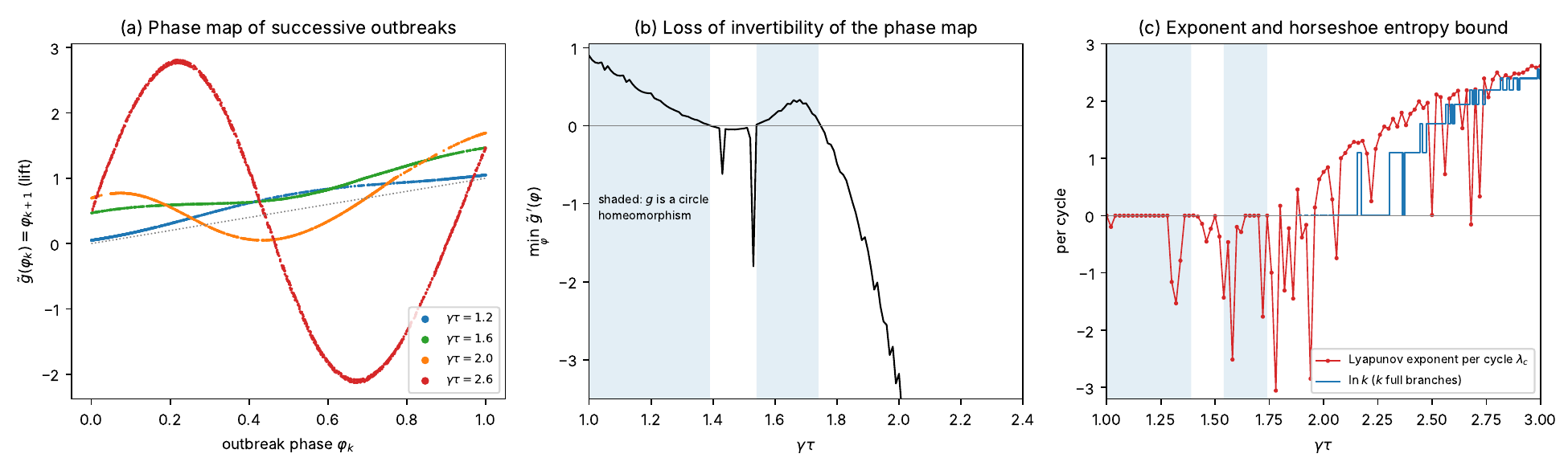}
\caption{Reduction to the phase map along the round length ($\rho=25$, $\nu=0.005$, $p=0.05$). (a) Lift of the phase map, $\varphi_{k+1}$ against the phase $\varphi_k$ of the previous outbreak (shifted by an integer), from the transients of 300 initial conditions: increasing for $\gamma\tau=1.2$ and $1.6$ (circle homeomorphism), folded for $\gamma\tau=2.0$ and $2.6$. (b) Minimum slope of the lift; shaded: $g$ is invertible. (c) Lyapunov exponent per cycle (red) and lower bound $\ln k$ for the topological entropy per cycle given by $k$ full branches of $g$ (blue); shaded as in (b).}
\label{fig:circlemap}
\end{figure}

Along $\rho$ at $\gamma\tau=3$, chaos persists above $\rNS=32.58$ in the bistable band and ends abruptly at $\rho\approx34.05$ with $\lambda_1=0.032$ per step, again a boundary crisis (Supplementary Material, Fig.~\ref{Sfig:route_rho}).

\begin{figure}[t]
\centering
\includegraphics[width=0.72\linewidth]{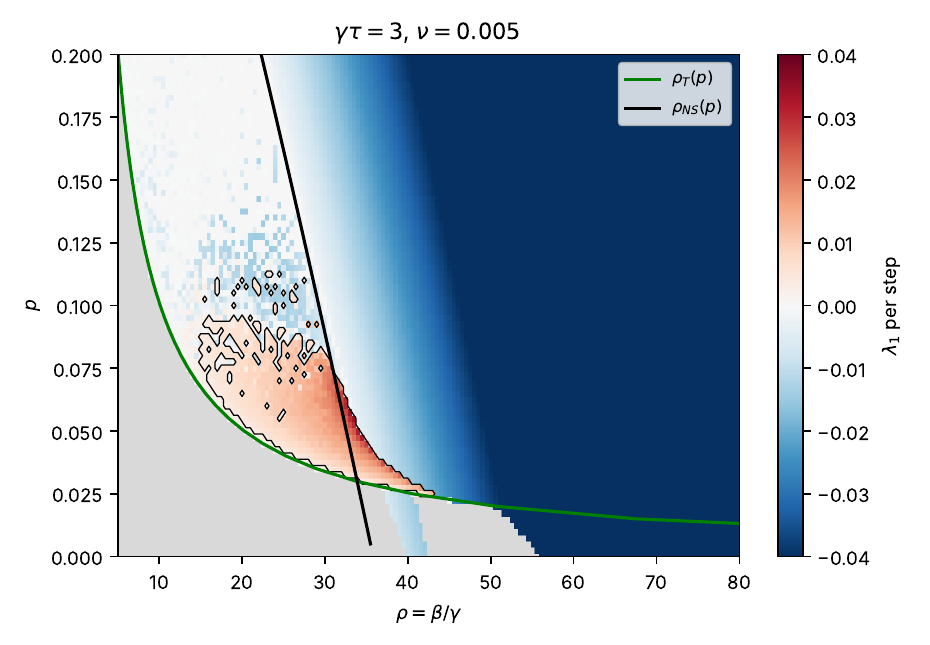}
\caption{Largest Lyapunov exponent per step in the plane $(\rho,p)$ for $\gamma\tau=3$ and $\nu=0.005$. Grey: extinction. Black contour: chaotic attractors ($\lambda_1>0.003$). Green: invasion threshold $\rT(p)$; black line: Neimark--Sacker curve $\rNS(p)$.}
\label{fig:rhop}
\end{figure}

Figure~\ref{fig:rhop} combines both parameters: $\lambda_1$ in the plane $(\rho,p)$ for $\gamma\tau=3$ and $\nu=0.005$. Chaos occurs for $0.025\le p\le0.113$ and $15\le\rho\le43$; about 5\% of the explored grid is chaotic, with a maximum exponent of $0.034$ per step. It is organised by the band $\rT(p)<\rho<\rNS(p)$, which opens only for $p\gtrsim0.03$ because $\rT(p)\propto1/p$ while $\rNS(p)$ varies slowly (from $35.5$ at $p=0.005$ to $23.7$ at $p=0.2$). Each horizontal line of Fig.~\ref{fig:rhop} is a sweep in $p$ of the kind studied in Section~\ref{sec:p}. The chaotic region grows as turnover slows down (Table~\ref{tab:ranges}): for $\nu=0.0005$ it covers $29\%$ of the plane $(\rho,\gamma\tau)\in[15,100]\times[0.5,60]$ and reaches the edges of that window, with exponents up to $0.066$ per step, and for $\nu=0.001$ and $\gamma\tau=10$ it extends to $0.02\le p\le0.13$, with exponents up to $0.094$ per step.

\begin{table}[t]
\centering\small
\begin{tabular}{@{}lll@{}}
\toprule
Parameters held fixed & Chaotic range found & Maximum $\lambda_1$ (per step)\\
\midrule
$\rho=46$, $\nu=0.002$, $\gamma\tau=5.5$ & $0.023\le p\le0.0555$ & $0.046$ ($p=0.054$)\\
$\gamma\tau=3$, $\nu=0.005$ & $0.025\le p\le0.113$, \ $15\le\rho\le43$ & $0.034$\\
$p=0.05$, $\nu=0.005$ & $1.7\lesssim\gamma\tau\lesssim7.6$, \ $20.8\le\rho\le46$ & $0.036$ ($\gamma\tau=5.3$, $\rho=25.8$)\\
$p=0.05$, $\nu=0.002$ & $2.6\lesssim\gamma\tau\le20^*$, \ $20.5\le\rho\le50^*$ & $0.048$ ($\gamma\tau=5.6$, $\rho=46$)\\
$p=0.05$, $\nu=0.0005$ & $3.5\lesssim\gamma\tau\le60^*$, \ $21\le\rho\le100^*$ & $0.066$ ($\gamma\tau=17.8$, $\rho=94$)\\
$\gamma\tau=10$, $\nu=0.001$ & $0.020\le p\le0.13$, \ $15\le\rho\le80$ & $0.094$ ($p=0.08$, $\rho=80$)\\
$p=0.05$, $\nu\ge0.02$ & none & ---\\
$p=0$, $p=0.5$ or $p=1$ & none & ---\\
permanent adoption ($\gamma=0$) & none & ---\\
\bottomrule
\end{tabular}
\caption{Parameter ranges with chaotic attractors ($\lambda_1>0.003$ per step) in the explored windows. The controls (last three rows) were scanned over $(\rho,\gamma\tau)\in[1,100]\times[0.05,20]$ and, for permanent adoption, over $\beta\in[0.005,2]$, $\tau\in[0.1,100]$ with $\mu=0.02$ and four values of $p$. An asterisk marks a range limited by the edge of the explored window.}
\label{tab:ranges}
\end{table}

\subsection{The two-stage parameter $p$ as bifurcation parameter}\label{sec:p}

\begin{figure}[t]
\centering
\includegraphics[width=\linewidth]{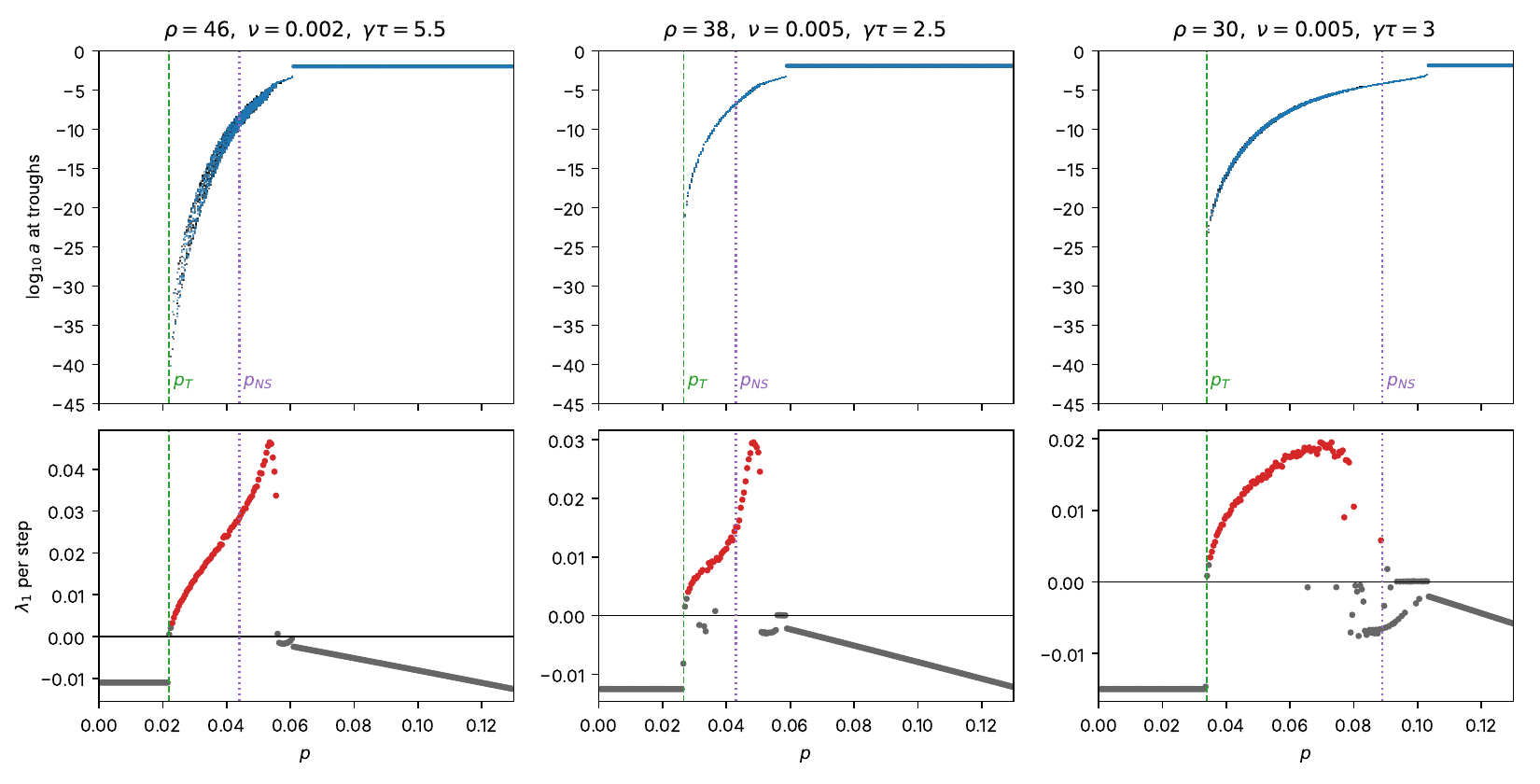}
\caption{The two-stage parameter $p$ as bifurcation parameter, with rates and round length fixed. Top: values of $\log_{10}a$ at the troughs of the attractor (black and blue: two initial conditions). Bottom: largest Lyapunov exponent per step (red: chaotic). Dashed green: invasion threshold $\pT$, Eq.~\eqref{eq:pT}; dotted purple: subcritical Neimark--Sacker point $\pNS$. For larger $p$, up to $p=1$, the equilibrium is stable.}
\label{fig:sweepp}
\end{figure}

Once the rounds are long enough, the two-stage parameter $p$ decides whether the oscillation becomes chaotic. It is the parameter that distinguishes two-stage from one-stage contagion, and the natural bifurcation parameter of the model. Figure~\ref{fig:sweepp} and Table~\ref{tab:sweepp} show the attractors of $F_\tau$ as $p$ varies from $0$ to $1$ with the rates and the round length fixed, for three parameter sets; the first, with $p=0.05$, is the reference point used below. In all three the sequence is the same.

\begin{enumerate}
\item For $p<\pT$, given by Eq.~\eqref{eq:pT}, adoption goes extinct after the first outbreak: the adoption-free state is stable ($r_0<1$) and outbreaks cannot re-ignite once the primed reservoir is depleted.
\item Just above $\pT$ the attractor is chaotic. The exponent grows from zero as $p$ increases, while the deepest troughs rise from $a\sim10^{-40}$ near $\pT$ to $a\sim10^{-5}$ at the upper end of the window.
\item The endemic equilibrium becomes stable at a subcritical NS point $\pNS$. Chaos persists slightly beyond $\pNS$, coexisting with the stable equilibrium, and disappears abruptly while $\lambda_1$ is still large. This is a boundary crisis~\cite{Grebogi1983}, in which the chaotic attractor collides with the boundary of the basin of the equilibrium.
\item For larger $p$, up to the one-stage limit $p=1$, all orbits computed converge to the equilibrium.
\end{enumerate}

\begin{table}[t]
\centering\small
\setlength{\tabcolsep}{4pt}
\begin{tabular}{@{}lcccccc@{}}
\toprule
fixed parameters & $\pT$ & $\pNS$ & chaotic window & chaotic & max.\ $\lambda_1$ & equilibrium\\
 & & & in $p$ & fraction & & for\\
\midrule
$\rho=46$, $\nu=0.002$, $\gamma\tau=5.5$ & 0.0219 & 0.0440 & $[0.023,\,0.0555]$ & 100\% & 0.046 & $p\ge0.061$\\
$\rho=38$, $\nu=0.005$, $\gamma\tau=2.5$ & 0.0266 & 0.0430 & $[0.028,\,0.0505]$ & 91\% & 0.030 & $p\ge0.059$\\
$\rho=30$, $\nu=0.005$, $\gamma\tau=3$ & 0.0338 & 0.0890 & $[0.035,\,0.0885]$ & 81\% & 0.020 & $p\ge0.104$\\
\bottomrule
\end{tabular}
\caption{Chaos as a function of $p$ (480 values of $p\in(0,1]$, two initial conditions). The chaotic fraction is that of the $p$ values inside the window; the remainder are periodic windows. Maximum exponents are per step.}
\label{tab:sweepp}
\end{table}

The window is narrow and lies at small $p$: chaos requires that most first contacts lead to the intermediate stage, but that a small fraction of single-contact adoptions remains, enough for outbreaks to re-ignite. For $p=0$ re-ignition is impossible (Section~\ref{sec:conditions}); for larger $p$ the two-stage mechanism loses weight and the dynamics approaches that of one-stage contagion, whose equilibrium is stable. Near $\pT$ the troughs become very deep, down to $a\sim10^{-40}$, and the intervals between outbreaks become long; as $p$ increases, the troughs become shallower and the outbreaks more frequent.

\label{sec:unique}To test whether several chaotic attractors coexist, we iterated 300 initial conditions drawn uniformly from the simplex at four points of the chaotic windows of Table~\ref{tab:sweepp}, and computed the full Lyapunov spectrum on the resulting attractor (Table~\ref{tab:unique}). At three of the points all 300 orbits converged to the same chaotic attractor: they have the same exponent to within $\pm0.001$, the same number of outbreaks and the same depth of troughs. At the fourth, close to a periodic window, 13\% of the orbits converged to a stable periodic orbit of period 98 (two outbreaks per period) and the rest to the chaotic attractor.

\begin{table}[t]
\centering\small
\begin{tabular}{@{}lccc@{}}
\toprule
point & final states of 300 orbits & $(\lambda_1,\lambda_2,\lambda_3)$ per step & $D_{KY}$\\
\midrule
$\rho=46$, $\nu=0.002$, $\gamma\tau=5.5$, $p=0.035$ & 300 chaotic & $(0.019,-0.027,-0.092)$ & 1.70\\
$\rho=46$, $\nu=0.002$, $\gamma\tau=5.5$, $p=0.05$ & 300 chaotic & $(0.037,-0.046,-0.256)$ & 1.81\\
$\rho=30$, $\nu=0.005$, $\gamma\tau=3$, $p=0.06$ & 300 chaotic & $(0.018,-0.034,-0.198)$ & 1.52\\
$\rho=30$, $\nu=0.005$, $\gamma\tau=3$, $p=0.08$ & 261 chaotic, 39 period 98 & $(0.010,-0.026,-0.404)$ & 1.39\\
\bottomrule
\end{tabular}
\caption{Uniqueness of the chaotic attractor: 300 random initial conditions per point.}
\label{tab:unique}
\end{table}

In all our tests, the chaotic attractor was therefore unique up to coexistence with the stable equilibrium (between $\pNS$ and the crisis) and, near periodic windows, with stable periodic orbits. It has a single positive exponent, so there is no hyperchaos, and a dimension between one and two: it is a thin folded band. Its geometry is the same throughout the window, a fast outbreak followed by a long recharge in which naive susceptibles are replenished by turnover; what changes with $p$ is the scale, with shallower troughs, shorter intervals between outbreaks and a lower dimension as $p$ increases. These tests sample the parameter window and the phase space, but they do not prove global uniqueness.

In summary, the attractor is a strange attractor of a three-dimensional dissipative map, with one positive Lyapunov exponent and a fractal dimension between one and two, a thin band stretched along the recharge phase and folded at each outbreak, as in a horseshoe. It is born from the break-up of an invariant closed curve, signalled by the loss of invertibility of the phase map of successive outbreaks (Section~\ref{sec:roundlength}), and destroyed in a boundary crisis. Because periodic windows are embedded in the chaotic range and stable periodic orbits can coexist with it, it is not hyperbolic; it is of the type called quasi-attractor by Afraimovich and Shilnikov~\cite{AfraimovichShilnikov1983}, typical of torus breakdown.

%======================================================================
\subsection{Time evolution and indicators of chaos}\label{sec:timeevol}

\begin{figure}[t]
\centering
\includegraphics[width=0.92\linewidth]{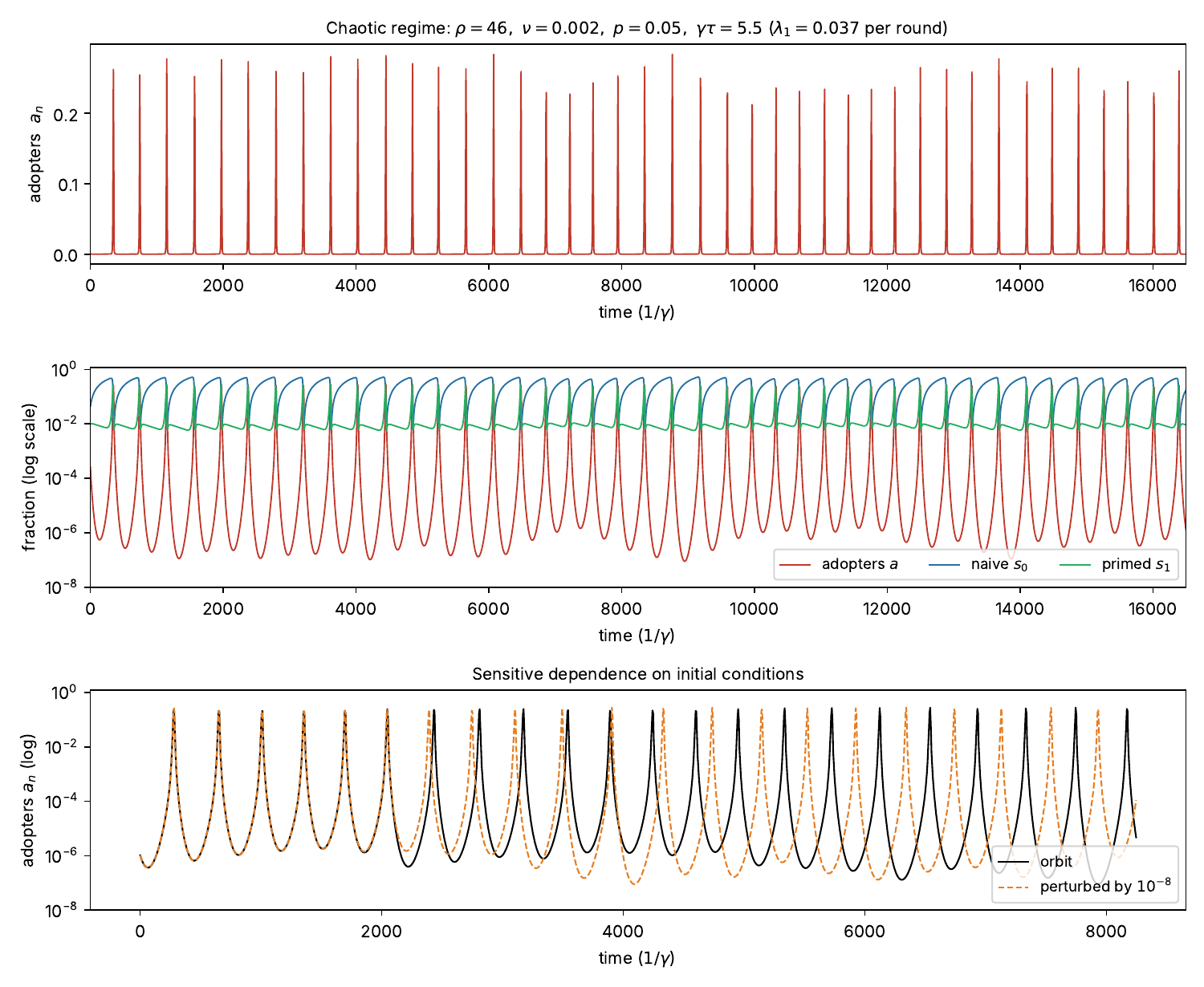}
\caption{Time evolution on the chaotic attractor at the reference point ($\rho=46$, $\nu=0.002$, $p=0.05$, $\gamma\tau=5.5$). Top: adopter fraction $a_n$ over 3000 rounds. Middle: $a_n$, $s_{0,n}$ and $s_{1,n}$ in logarithmic scale. Bottom: two orbits started $10^{-8}$ apart.}
\label{fig:timeseries}
\end{figure}

Figure~\ref{fig:timeseries} shows the time evolution on the attractor at the reference point, $\rho=46$, $\nu=0.002$, $p=0.05$, $\gamma\tau=5.5$. Adoption proceeds by irregular outbreaks, with peaks of $a\approx0.21$--$0.29$, separated by long recharge phases in which adopters fall to $a\sim10^{-6}$--$10^{-7}$ while naive susceptibles are replenished by turnover; the minimum adopter fraction on the attractor is about $7\times10^{-8}$. Neither the heights of the peaks nor the intervals between outbreaks repeat. Two orbits that start $10^{-8}$ apart remain indistinguishable for about six outbreaks and then become completely different (Fig.~\ref{fig:timeseries}, bottom). The Lyapunov spectrum, computed by QR re-orthonormalisation over $2\times10^5$ steps, is
\[
(\lambda_1,\lambda_2,\lambda_3)=(0.0373,\,-0.0462,\,-0.2564)\ \text{per step},\qquad D_{KY}=1+\lambda_1/|\lambda_2|\approx1.81 .
\]
Although $\lambda_1$ is small per round, one outbreak--recharge cycle lasts about 70 rounds, over which perturbations grow by a factor $e^{0.037\cdot70}\approx13$. Since $\pNS=0.044$ for these rates, the reference point lies in the band in which the chaotic attractor coexists with the stable equilibrium; the basin of the equilibrium is small, and none of 300 random initial conditions reached it.

\begin{figure}[t]
\centering
\includegraphics[width=\linewidth]{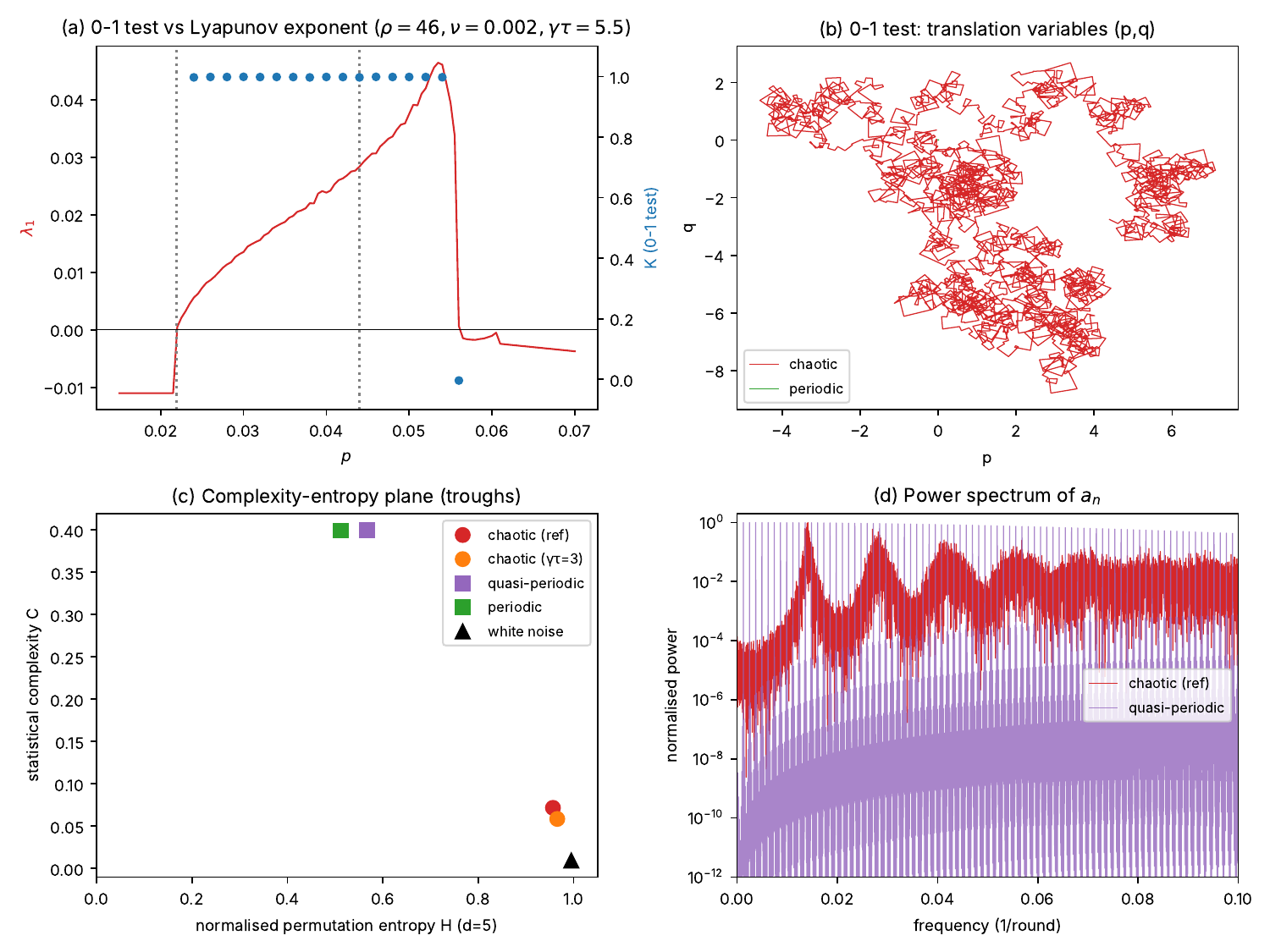}
\caption{Indicators of chaos. (a) 0-1 test $K$ (blue) and largest Lyapunov exponent per round (red) along $p$ at the reference rates; dotted lines: $\pT$ and $\pNS$. (b) Translation variables of the 0-1 test: diffusive for the chaotic orbit (red); the periodic orbit (green) stays bounded near the origin. (c) Complexity--entropy plane of the trough series. (d) Power spectra of $a_n$ for the chaotic and the quasi-periodic orbit.}
\label{fig:indicators}
\end{figure}

\label{sec:indicators}A positive Lyapunov exponent is the standard signature of chaos, but its estimate can be affected by numerical issues and it does not by itself exclude long quasi-periodic or periodic transients. We therefore applied several independent diagnostics (Fig.~\ref{fig:indicators}, Table~\ref{tab:indicators}). Except for the spectrum, they were computed from the sequence of trough values $\log_{10}a$, one value per cycle, which removes the trivial oscillation within each cycle; we compare chaotic points with a quasi-periodic point ($\rho=25$, $\nu=0.005$, $p=0.05$, $\gamma\tau=1.2$), a periodic point in a periodic window ($\rho=30$, $\nu=0.005$, $p=0.079$, $\gamma\tau=3$) and white noise.

\begin{itemize}
\item \emph{0-1 test.} The test of Gottwald and Melbourne~\cite{Gottwald2004,Gottwald2009} drives a two-dimensional system with the observable and measures whether it spreads diffusively ($K\approx1$, chaos) or stays bounded ($K\approx0$, regular dynamics). It gives $K=0.998$ at the reference point and $K\approx0$ for the quasi-periodic and periodic points. Along $p$ at the reference rates, $K\approx1$ for $0.024\le p\le0.054$ and $K\approx0$ at $p=0.056$, reproducing the chaotic window obtained from $\lambda_1$ to within the resolution of the grid (Fig.~\ref{fig:indicators}a,b).
\item \emph{Power spectrum.} The spectrum of $a_n$ is broadband at the reference point, with only $20\%$ of the power in the $1\%$ most intense frequencies, whereas the quasi-periodic orbit concentrates $97\%$ of its power in discrete peaks (Fig.~\ref{fig:indicators}d).
\item \emph{Ordinal patterns.} The normalised permutation entropy $H$ and the statistical complexity $C$ of Bandt and Pompe~\cite{BandtPompe2002,Rosso2007} (patterns of length 5) place the chaotic points far from regular dynamics ($H\approx0.5$, $C\approx0.4$) and, through a complexity eight times that of white noise, away from noise (Fig.~\ref{fig:indicators}c). The entropy is close to that of noise because one value per cycle is sampled after perturbations have been amplified by a factor of about 13.
\item \emph{Determinism.} A nearest-neighbour prediction of the next trough from the previous three has a normalised error of $1.09$, against $1.41\pm0.03$ for shuffled surrogates of the same series: the sequence is deterministic, although a single observable per cycle does not capture the full three-dimensional state.
\item \emph{Geometry.} The correlation dimension, $D_2=1.75$--$1.85$ depending on the fitting window, agrees with $D_{KY}=1.81$. Newton's method found saddle periodic orbits with one outbreak for every period from 58 to 79 steps and with two outbreaks for periods 120--129, each with one unstable and two stable multipliers and growth rates $\ln|\Lambda_u|/k=0.03$--$0.05$ per step, comparable with $\lambda_1$. The unstable manifold of the period-70 orbit is stretched and folded as in a horseshoe, and the inter-outbreak intervals form a symbolic dynamics with positive topological entropy, about $0.031$ per step (Supplementary Material, Sections~\ref{Ssec:upo}--\ref{Ssec:d2}).
\item \emph{Numerical robustness.} In double and extended (80-bit) precision, over $4\times10^5$ steps, with two initial conditions and perturbations of $10^{-9}$ and $10^{-7}$, $\lambda_1=0.0368$--$0.0376$ at the reference point, $0.0080$--$0.0083$ at a weakly chaotic point ($\rho=25$, $\nu=0.005$, $p=0.05$, $\gamma\tau=3$) and $0.0309$--$0.0317$ near the crisis ($\rho=34$). Because the troughs are deep, we also computed $\lambda_1$ in the coordinates $(s_0,s_1,\ln a)$, in which all variables are of order one; the result, $0.0381$, agrees with the value in the original coordinates, as it must for an invariant quantity. Small values of $a$ do not reduce the relative precision of floating-point arithmetic, and the map involves no cancellation in the equation for $a$.
\end{itemize}

\begin{table}[t]
\centering\small
\begin{tabular}{@{}lccc@{}}
\toprule
series of troughs & 0-1 test $K$ & entropy $H$ & complexity $C$\\
\midrule
chaotic, reference point & 0.998 & 0.959 & 0.068\\
chaotic, $\rho=30$, $\nu=0.005$, $p=0.06$, $\gamma\tau=3$ & 0.997 & 0.967 & 0.055\\
quasi-periodic, $\rho=25$, $\nu=0.005$, $p=0.05$, $\gamma\tau=1.2$ & 0.004 & 0.566 & 0.399\\
periodic, $\rho=30$, $\nu=0.005$, $p=0.079$, $\gamma\tau=3$ & $-0.004$ & 0.512 & 0.399\\
white noise & 0.998 & 0.996 & 0.008\\
\bottomrule
\end{tabular}
\caption{Indicators of chaos computed from the sequence of trough values (about 2500--3000 values per series). $H$ is the normalised permutation entropy and $C$ the statistical complexity for ordinal patterns of length 5.}
\label{tab:indicators}
\end{table}

\subsection{Short rounds, and what Katriel's model does in the chaotic window}\label{sec:katrielwindow}\label{sec:why}

We now turn from what the map does to why it does it. By Proposition~\ref{prop:consistency}, $F_\tau$ behaves for small $\tau$ as a smooth one-step method of order one for Eq.~\eqref{eq:ode}. Classical results on one-step discretisations then constrain its dynamics for short rounds. If the ODE has a hyperbolic attracting periodic orbit $\Gamma$, then for all sufficiently small $\tau$ the map has an attracting invariant closed curve $\Gamma_\tau$ that converges to $\Gamma$ as $\tau\to0$ and inherits its stability~\cite{Beyn1987}. More generally, uniformly asymptotically stable attracting sets of the ODE have nearby attracting sets of the discretisation~\cite{KloedenLorenz1986}. As in Section~\ref{sec:roundlength}, the dynamics on an invariant closed curve is that of a circle homeomorphism, which has zero topological entropy; its Lyapunov exponent along the curve is zero (quasi-periodic motion) or negative (periodic, mode-locked motion). Hence, near the hyperbolic attractors of Katriel's model, $F_\tau$ cannot have chaotic attractors when the rounds are short. The bottom of Fig.~\ref{fig:classes} and the invertible range of the phase map in Fig.~\ref{fig:circlemap}b are the numerical counterpart of this statement.

These results explain why there is no chaos for short rounds, but not what Katriel's model does where the map is chaotic. To see this, we integrated Eq.~\eqref{eq:ode} with the same rates as in Table~\ref{tab:sweepp}, over the same range of $p$ (Table~\ref{tab:katriel}). The invasion threshold of the ODE, $p=(1+\nu)/\rho$, coincides with $\pT$ to within $1\%$. Above it, and throughout the chaotic window of the map, the ODE has a periodic relaxation cycle; the equilibrium becomes stable only at a Hopf point $p_H\approx0.31$--$0.41$, far beyond the end of the chaotic window.

\begin{table}[t]
\centering\small
\begin{tabular}{@{}lcccc@{}}
\toprule
& \multicolumn{2}{c}{Katriel's ODE} & \multicolumn{2}{c}{map $F_\tau$}\\
\cmidrule(lr){2-3}\cmidrule(l){4-5}
fixed parameters & $p_H$ & troughs (periodic cycle) & troughs (chaos) & oscillation ends\\
\midrule
$\rho=46$, $\nu=0.002$, $\gamma\tau=5.5$ & 0.41 & $10^{-164}$--$10^{-44}$ & $10^{-39}$--$10^{-4.8}$ & $p\approx0.06$\\
$\rho=38$, $\nu=0.005$, $\gamma\tau=2.5$ & 0.31 & $10^{-57}$--$10^{-21}$ & $10^{-19}$--$10^{-4.4}$ & $p\approx0.06$\\
$\rho=30$, $\nu=0.005$, $\gamma\tau=3$ & 0.35 & $10^{-63}$--$10^{-16}$ & $10^{-22}$--$10^{-4.2}$ & $p\approx0.10$\\
\bottomrule
\end{tabular}
\caption{Katriel's ODE and the map with the same rates, over the chaotic window of the map. In the whole window the ODE is on a periodic relaxation cycle, and it oscillates up to its Hopf point $p_H$. Troughs are given at the lower and upper ends of the chaotic window of Table~\ref{tab:sweepp}. At the reference rates the period of the ODE cycle decreases from $8450/\gamma$ at $p=0.023$ to $503/\gamma$ at $p=0.0555$.}
\label{tab:katriel}
\end{table}

In the window where the map is chaotic, Katriel's ODE therefore has a periodic relaxation cycle with troughs of $10^{-44}$--$10^{-164}$ at the reference rates. The round structure does not create the oscillation: it shortens the recharge phase, raising the troughs by tens of orders of magnitude, turns the cycle into a chaotic attractor, and shrinks the oscillatory range in $p$ from $(\pT,p_H\approx0.41)$ to $(\pT,\approx0.056)$. Most of Katriel's oscillatory range, $0.06\lesssim p<0.41$, becomes a stable equilibrium when updating is synchronous; chaos occupies its small-$p$ end, where the two-stage mechanism dominates.

\subsection{Mechanism and conditions for chaos}

\begin{figure}[t]
\centering
\includegraphics[width=\linewidth]{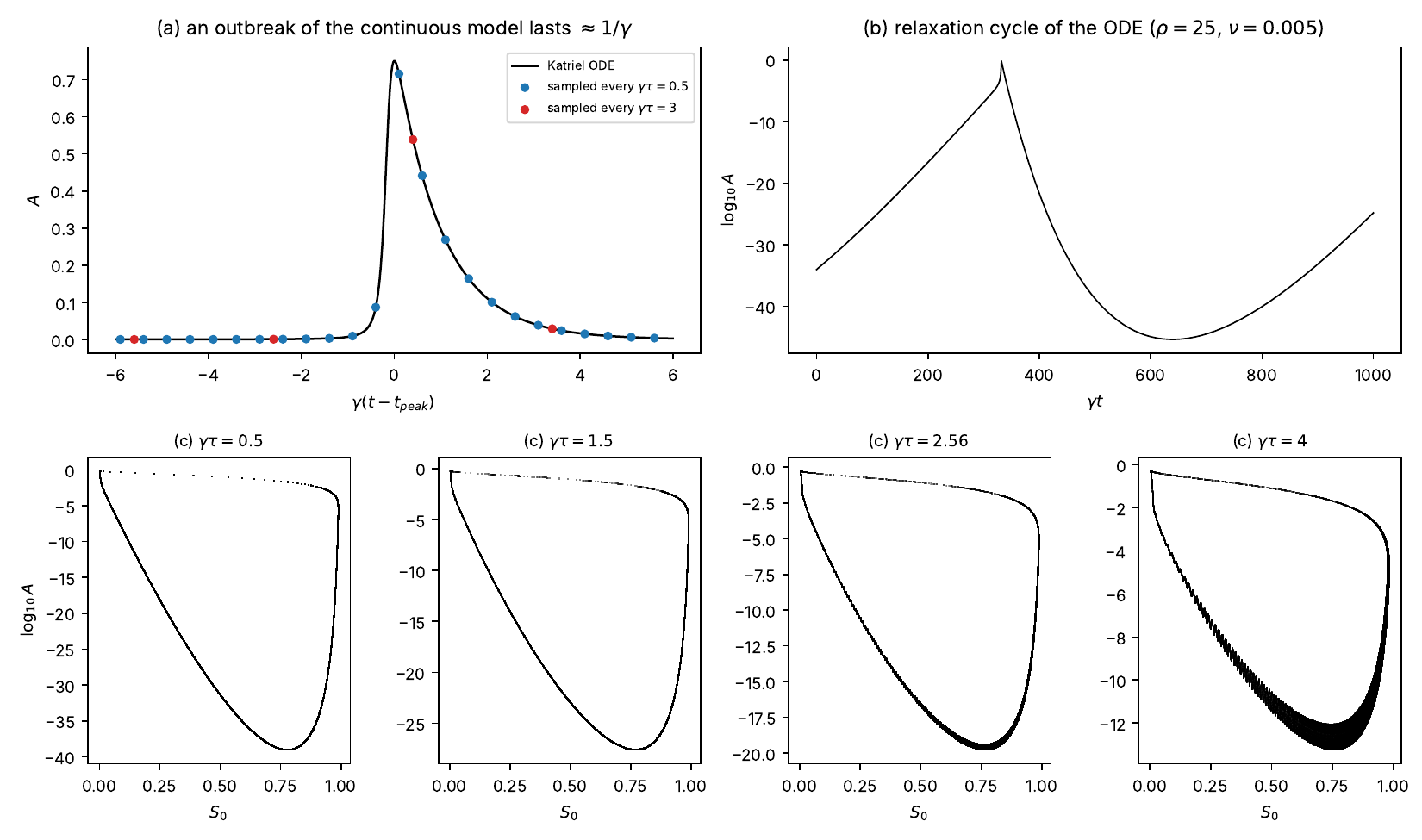}
\caption{Mechanism ($\rho=25$, $\nu=0.005$, $p=0.05$). (a) An outbreak of Katriel's ODE in linear scale, sampled with rounds $\gamma\tau=0.5$ (blue) and $3$ (red). (b) The relaxation cycle of the ODE in logarithmic scale. (c) Attractors of $F_\tau$ for increasing $\gamma\tau$: closed curves for $\gamma\tau=0.5$ and $1.5$; a chaotic attractor with folded layers for $\gamma\tau=2.56$ and $4$.}
\label{fig:mechanism}
\end{figure}

The mechanism is a mismatch between two time scales: the round length and the duration of an outbreak. The relaxation cycle of Eq.~\eqref{eq:ode} has two very different time scales. Recharge of naive susceptibles by turnover takes hundreds of adoption times; in our examples the period is $180/\gamma$ to $8500/\gamma$, and $a$ falls to $10^{-16}$--$10^{-164}$. The outbreak itself is fast. Its full width at half maximum is $0.95/\gamma$--$1.05/\gamma$ in all cases we computed ($\rho=22$--$46$, $\nu=0.002$--$0.005$), essentially the mean adoption time (Fig.~\ref{fig:mechanism}a,b).

As long as $\gamma\tau\ll1$, the map follows the outbreak with many steps and its attractor is a smooth invariant curve. When $\gamma\tau$ approaches $2$, i.e.\ about twice the outbreak width, an outbreak is covered by at most one or two rounds (red points in Fig.~\ref{fig:mechanism}a). The size of the realised outbreak and the state it leaves behind then depend sharply on where the threshold crossing falls relative to the sequence of rounds, and that phase changes from one cycle to the next because the recharge time is not a multiple of $\tau$. In addition, adoptions made during a round act only from the next one, so that contagion during an outbreak is concentrated in a few successive generations of adopters. The phase map of successive outbreaks becomes non-invertible (Fig.~\ref{fig:circlemap}) and the invariant curve wrinkles and breaks up (Fig.~\ref{fig:mechanism}c, $\gamma\tau=2.56$ and $4$). The same synchronisation stabilises the endemic equilibrium, which is why the NS curve descends as $\tau$ grows and the oscillatory range in $p$ shrinks. The onset of chaos is therefore set by the ratio of the round length to the outbreak duration, which is close to one adoption time in all cases, and not by the period of the cycle, which is hundreds to thousands of times longer: this is why chaos appears at a similar round length, $\gamma\tau\approx1.7$--$2.6$, for both turnover ratios in Fig.~\ref{fig:classes}, although their recharge times differ by a factor of about $2.5$.

\label{sec:conditions}Taken together, these results show that chaos requires four ingredients. (a) The two-stage mechanism with turnover, which produces the relaxation oscillation, as in Katriel's model. (b) A small but positive probability $p$ of single-contact adoption. Near $a=0$, Eq.~\eqref{eq:map} gives $a_{n+1}\approx r_na_n$ with
\[
r_n=\big(1-e^{-(1+\nu)\tau}\big)\,\frac{\rho}{1+\nu}\,\frac{1-e^{-\nu\tau}}{\nu\tau}\,\big(ps_{0,n}+s_{1,n}\big)+e^{-(1+\nu)\tau}.
\]
During the recharge phase $s_1$ decays, because priming requires adopters, while $s_0$ tends to one; hence $r_n$ tends to the eigenvalue $r_0$ of Proposition~\ref{prop:threshold}, and re-ignition from a deep trough is possible if and only if $p>\pT$. For $p=0$, $r_0=e^{-(1+\nu)\tau}<1$ for every $\rho$; orbits either go extinct after one outbreak or, for $\rho\gtrsim40$ at $\gamma\tau=3$, converge to a stable endemic equilibrium, and no sustained oscillation or chaos was found. If $p$ is too large, the two-stage mechanism loses weight and the equilibrium is stable. (c) Rounds of at least about twice the adoption time, which synchronise generations and break the oscillation into chaos. (d) Small turnover per round, $\mu\tau\lesssim0.04$; otherwise the equilibrium becomes stable.

\subsection{Robustness of the results}\label{sec:notspurious}

The evidence for chaos is numerical, but it is varied and consistent (Section~\ref{sec:indicators}). The positive exponent is insensitive to precision, coordinates, run length, initial condition and size of the perturbation; the 0-1 test, the spectrum, ordinal-pattern complexity and a determinism test against surrogates detect the same chaotic window independently of the exponent; the correlation dimension agrees with the Kaplan--Yorke dimension; and saddle periodic orbits of every period in the observed range have growth rates comparable with $\lambda_1$. Recomputing the $(\rho,\gamma\tau)$ plane with eight random initial conditions and longer runs ($4\times10^4+4\times10^4$ steps) reproduced all 168 chaotic points of a subgrid and added two, and the main results are robust to the choice of exposure closure (P5) (Remark~\ref{rem:profile}). Once the reduction to the phase map of Section~\ref{sec:roundlength} is accepted, part of the argument is rigorous: a circle homeomorphism has zero topological entropy, so there is no chaos while $g$ is invertible, and $k$ disjoint arcs mapped across the circle force a topological entropy of at least $\ln k$. What remains numerical is the reduction itself and the covering of the arcs. Both are of the kind that can be validated with interval arithmetic by covering relations for h-sets~\cite{Zgliczynski2004}: the coverings have margins of $0.05$--$0.57$ in units of the circle, at least five times the scatter of the reduction and more than eight orders of magnitude above the floating-point error (at most $2\times10^{-10}$ in the phase over three cycles, by comparison with 80-bit arithmetic), and the strong transverse contraction makes the stable directions easy to control. A computer-assisted proof of the horseshoe therefore appears feasible; carrying it out is left for future work.

The chaos is also specific to the two-stage mechanism and to the round structure. The map is the large-population limit of a process derived from (P1)--(P5), not a numerical scheme for Eq.~\eqref{eq:ode}, and it preserves positivity, the population and the invasion threshold for every $\tau$. For the same rounds, one-stage contagion ($p=1$), $p=0.5$, $p=0$ and permanent adoption are never chaotic (Table~\ref{tab:ranges}), and for fixed rates chaos exists only in a window of small $p$. The route to chaos, break-up of an invariant curve followed by a boundary crisis, differs from the period-doubling route typical of discrete-time epidemic maps~\cite{Allen1994,Hu2014}.

%======================================================================
\section{Discussion}

We have shown that round-based updating can turn a regular oscillation into chaos. Starting from elementary processes and five explicit principles, we derived a discrete-time two-stage contagion model whose deterministic large-population map reduces to Katriel's equations for short rounds, and showed that it is chaotic in a well-defined region of parameter space. The central message is that a continuous-time relaxation oscillator becomes chaotic when updates are forced to occur in synchronous rounds whose duration is comparable to the outbreak duration. The essential additional ingredient is that the population is updated in rounds: each individual makes at most one transition per round, and adoptions take effect from the next round. The within-round exposure closure (P5) is introduced to obtain closed-form transition probabilities, and the results are robust to this approximation. Round-based updating is characteristic of many social processes, among them elections, iterated games, periodic meetings and seasonal or release cycles, and the round length is then a parameter of the process, not of a numerical method.

The central parameter of the model is $p$, the probability of adoption after a single contact, which sets how much of the contagion goes through the intermediate stage. With all rates and the round length fixed, $p$ organises the dynamics: extinction below the invasion threshold $\pT$, chaos in a narrow window above it, and a stable equilibrium beyond a subcritical Neimark--Sacker point and a boundary crisis, up to the one-stage limit $p=1$. Chaos is thus a property of two-stage contagion with a small residual probability of direct adoption, and disappears in the one-stage limit.

In all our tests, within the chaotic window the attractor was unique up to coexistence with the stable equilibrium (between $\pNS$ and the crisis) and, near periodic windows, with stable periodic orbits; it has a single positive exponent and the same outbreak--recharge geometry throughout. Besides the Lyapunov exponent, chaos is detected by the 0-1 test, a broadband spectrum, ordinal-pattern complexity, a determinism test and a correlation dimension that matches the Kaplan--Yorke dimension; it is robust to every numerical test we applied; and its periodic orbits, its unstable manifold and the phase map of successive outbreaks show a horseshoe. In the window where the map is chaotic, Katriel's ODE has a periodic relaxation cycle with astronomically deep troughs; the round structure does not create the oscillation but shortens its recharge phase, turns the cycle into a chaotic attractor, and shrinks the oscillatory range in $p$ from $(\pT,p_H\approx0.41)$ to $(\pT,\approx0.056)$. Theorems on one-step discretisations explain why chaos is absent for short rounds; the phase map of successive outbreaks shows that the invariant curve breaks down when this map loses invertibility, at $\gamma\tau\approx1.4$ along the line studied, and the numerical analysis locates the round length at which chaos appears: about twice the mean adoption time, the scale at which outbreaks are no longer resolved and successive generations of adopters become synchronised.

Compared with other discrete and multi-stage contagion models, which reach absorbing or stationary states, and with discrete epidemic maps, in which chaos arises by period doubling, the present model shows a different phenomenon: deterministic chaos produced by the interplay of two-stage contagion, renewal of susceptibles and synchronous updating, reached through the break-up of an invariant curve. Moreover, the troughs of the chaotic attractors are many orders of magnitude shallower than those of Katriel's relaxation cycle, which suggests that the chaotic regime is much less fragile to demographic noise. Because the mechanism requires only a slow--fast relaxation cycle and synchronous updating on the time scale of the fast phase, we expect it to operate in other round-based systems with relaxation dynamics; testing this is left for future work.

The results are numerical and limited to the explored windows; the parameter ranges in Tables~\ref{tab:sweepp}--\ref{tab:ranges} should be read in that sense. Natural continuations are a two-parameter continuation of the torus-breakdown and crisis bifurcations in the $(p,\tau)$ plane; a computer-assisted proof of the horseshoe of the phase map (Section~\ref{sec:roundlength}), for which the margins found here leave ample room; the study of the underlying stochastic chain in finite populations, where demographic noise acts in the deep troughs; empirical calibration of $\gamma\tau$, $\mu\tau$ and $p$ for specific processes organised in rounds; and extensions to networks and spatial settings~\cite{Meacci2025}.

\appendix
\section{Numerical methods}\label{app:num}
All computations use double precision, with $1-e^{-z}$ evaluated as \texttt{-expm1(-z)}; robustness checks used 80-bit extended precision. Lyapunov exponents on grids were obtained from a twin trajectory at distance $10^{-9}$ renormalised every step; spectra from QR re-orthonormalisation of finite-difference Jacobians. Grids: $(\rho,\gamma\tau)\in[10,50]\times[0.05,20]$ with $161\times100$ points (logarithmic in $\tau$); $(\rho,p)\in[5,80]\times[0,0.2]$ with $151\times81$ points at $\gamma\tau=3$ ($2\times10^4$ transient and averaging steps). Sweeps in $p$: 400 values with spacing $5\times10^{-4}$ in $(0,0.2]$ and 80 values with spacing $0.01$ in $(0.2,1]$, with $6\times10^4$ transient and $6\times10^4$ averaging steps. Initial conditions $(0.3,0.1,0.05)$ and $(0.6,0.3,0.2)$; the maximum exponent over both is reported. Uniqueness tests used 300 initial conditions drawn uniformly on the simplex, $4\times10^4$ transient and $6\times10^4$ averaging steps; periodic orbits were identified by return to within a relative distance $10^{-6}$. Endemic equilibria were computed from Proposition~\ref{prop:branch} or by bisection on the scalar equilibrium condition for $a$, and NS points by bisection on the spectral radius; Lyapunov coefficients from Kuznetsov's formula with finite-difference differentials. Saddle periodic orbits were computed by Newton's method on $F_\tau^k(x)=x$ from near-recurrences of a long orbit, and the correlation dimension by the Grassberger--Procaccia method with a Theiler window. Katriel's ODE was integrated with LSODA (relative tolerance $10^{-10}$) in the variables $(s_0,s_1,\ln a)$ to resolve deep troughs, and the Hopf point $p_H$ was located from the eigenvalues of the endemic equilibrium. Trough series were extracted as local minima of $\log_{10}a_n$. Outbreak phases were obtained from the upward crossings of $\log_{10}a=-3$, interpolating linearly between steps; derivatives of the phase map along orbits were computed by complex-step differentiation of $F_\tau$ (step $10^{-20}$), renormalised at each crossing; minimum slopes are $0.5\%$ quantiles over 300 transients of three cycles after two cycles of transient, and exponents per cycle were averaged over $1500$ cycles. The 0-1 test used the correlation method with $100$ random frequencies $c\in(\pi/5,4\pi/5)$ and the median of $K$; permutation entropy and statistical complexity used patterns of length 5; nearest-neighbour prediction used delay vectors of three troughs, the first half of the series as library and 20 shuffled surrogates. Exponents in the coordinates $(s_0,s_1,\ln a)$ were computed with a twin trajectory over $2\times10^5$ steps. Additional scans for $\nu=0.0005$ and $\nu=0.001$ used $86\times60$ and $74\times60$ grids with $2\times10^4$ transient and averaging steps.

\section*{Supplementary Material}
The Supplementary Material, at the end of this document (Sections~\ref{Ssec:route}--\ref{Ssec:d2}), contains the routes to chaos along $\rho$ and $\gamma\tau$, the reduction to the phase map with the details of the torus breakdown and of the one-dimensional horseshoe, and the saddle periodic orbits, unstable manifold, symbolic dynamics and correlation dimension at the reference point.

\section*{Data and code availability}
No empirical data were used. The Python code used for all computations and figures is available from the corresponding author upon reasonable request.

\section*{Author contributions (CRediT)}
\textbf{F.J.M., J.C.N.\ and L.M.} contributed equally to Conceptualization, Methodology, Software, Formal analysis, Investigation, Visualization, Writing -- original draft, and Writing -- review \& editing.

\section*{Acknowledgments}
L.M. acknowledges the support of the Gruppo Nazionale per la Fisica Matematica (GNFM) of the Istituto Nazionale di Alta Matematica (INdAM).

\section*{Funding}
This research received no specific grant from any funding agency in the public, commercial, or not-for-profit sectors.

\section*{Declaration of competing interest}
The authors declare that they have no known competing financial interests or personal relationships that could have appeared to influence the work reported in this paper.

\section*{Declaration of generative AI and AI-assisted technologies in the manuscript preparation process}
The research questions, hypotheses and modelling choices were formulated by the authors. During the work, the authors used Claude (Anthropic) as an assistant to test these hypotheses numerically, by writing and checking the Python code used for the computations and figures, to search the literature, and to help draft and edit the text. All figures were produced by this code from computed data; no generative AI tool was used to create or alter images. The authors reviewed and verified all code, results and text, and take full responsibility for the content of the publication.

\clearpage
\setcounter{section}{0}\setcounter{figure}{0}\setcounter{table}{0}\setcounter{equation}{0}
\renewcommand{\thesection}{S\arabic{section}}
\renewcommand{\thefigure}{S\arabic{figure}}
\renewcommand{\thetable}{S\arabic{table}}
\renewcommand{\theequation}{S\arabic{equation}}
\begin{center}{\Large\bfseries Supplementary Material}\end{center}
\addcontentsline{toc}{section}{Supplementary Material}

\noindent This Supplementary Material documents results that are summarised in the main text: the routes to chaos along the contact ratio and along the round length, including the reduction to a circle map (Section~\ref{Ssec:route}), and the periodic-orbit structure, unstable manifold, symbolic dynamics and correlation dimension of the chaotic attractor at the reference point (Sections~\ref{Ssec:upo}--\ref{Ssec:d2}). All results refer to the deterministic map $F_\tau$; time is measured in units of $1/\gamma$, $\rho=\beta/\gamma$, $\nu=\mu/\gamma$, and Lyapunov exponents are given per step (per round).
%======================================================================
\section{Routes to chaos}\label{Ssec:route}

\subsection{Along the contact ratio}

\begin{figure}[h]
\centering
\includegraphics[width=0.8\linewidth]{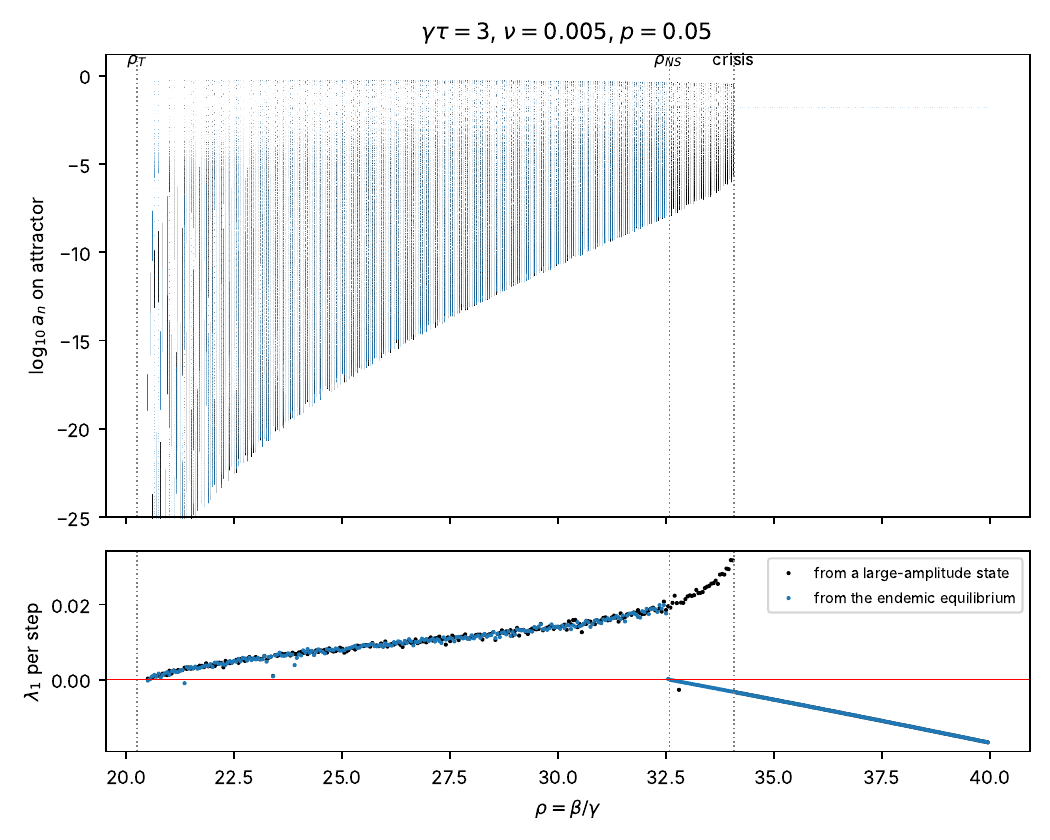}
\caption{Route along $\rho$ for $\gamma\tau=3$, $\nu=0.005$, $p=0.05$. Top: $\log_{10}a_n$ on the attractor (black: from a large-amplitude state; blue: from the endemic equilibrium). Bottom: largest Lyapunov exponent per step. Dashed lines: $\rT$, $\rNS$ and the crisis.}
\label{Sfig:route_rho}
\end{figure}

We fixed $\gamma\tau=3$, $\nu=0.005$, $p=0.05$ and varied $\rho\in[20.5,40]$ in steps of $0.05$, following the attractors from a large-amplitude state $(0.3,0.1,0.05)$ and from a small perturbation of the endemic equilibrium (Fig.~\ref{Sfig:route_rho}). The reference values are $\rT=20.25$ and $\rNS=32.58$ (subcritical).

Just above $\rT$ the oscillation is weakly chaotic, with $\lambda_1\lesssim0.003$, and $\lambda_1$ grows almost linearly with $\rho$, from $0.003$ near $\rho=21.4$ to $0.02$ near $\rho=32$; troughs deepen towards $\rT$, down to $a\sim10^{-25}$. Continued from the equilibrium, orbits remain there for $\rho\ge32.65$, in agreement with $\rNS$. Continued from the large-amplitude state, the chaotic attractor persists above $\rNS$ and coexists with the stable equilibrium for $32.6\lesssim\rho\lesssim34.05$, with an isolated periodic window (period 406) at $\rho=32.80$. At $\rho=34.05$ the chaotic attractor still has $\lambda_1=0.032$; at $\rho=34.10$ all orbits converge to the equilibrium ($\lambda_1=-0.0033$). The abrupt disappearance of an attractor with a large exponent, inside the bistable band, identifies a boundary crisis: the chaotic attractor collides with the boundary of the basin of the equilibrium, presumably the unstable invariant curve created at the subcritical Neimark--Sacker bifurcation.

\subsection{Along the round length}

\begin{figure}[h]
\centering
\includegraphics[width=0.8\linewidth]{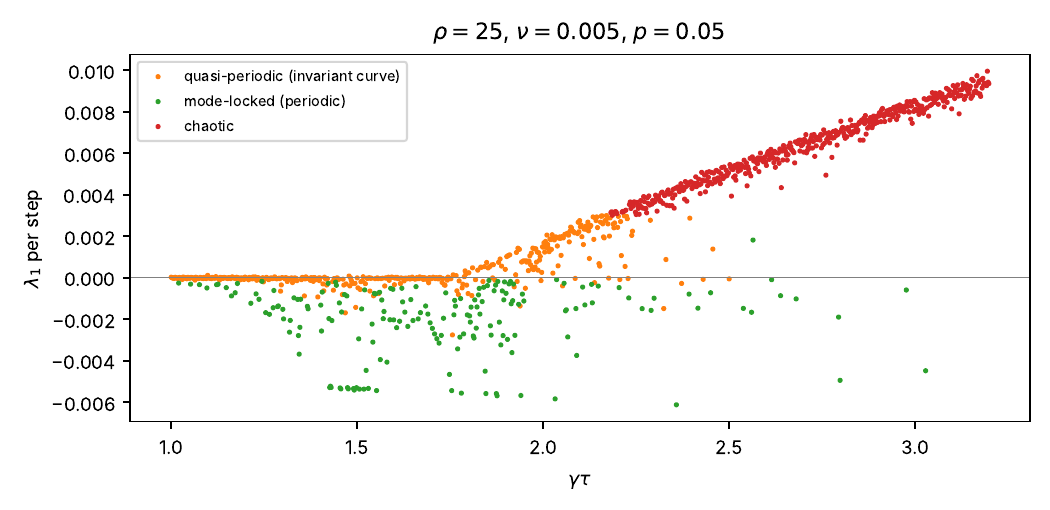}
\caption{Largest Lyapunov exponent per step against $\gamma\tau$ for $\rho=25$, $\nu=0.005$, $p=0.05$. Orange: quasi-periodic; green: mode-locked; red: chaotic.}
\label{Sfig:route_tau}
\end{figure}

We fixed $\rho=25$, $\nu=0.005$, $p=0.05$ and varied $\gamma\tau\in[1,3.2]$ in steps of $0.002$ (1100 values), with $3\times10^4$ transient and $6\times10^4$ averaging steps; orbits were tested for periodicity up to period 1500. Table~\ref{Stab:tau} gives the fraction of each type of attractor by interval of $\gamma\tau$, and Fig.~\ref{Sfig:route_tau} the exponent.

\begin{table}[h]
\centering\small
\begin{tabular}{@{}lcccc@{}}
\toprule
$\gamma\tau$ & quasi-periodic & mode-locked (periodic) & chaotic & steps per outbreak\\
\midrule
$[1.0,1.1)$ & 0.92 & 0.08 & 0 & 923\\
$[1.4,1.5)$ & 0.56 & 0.44 & 0 & 674\\
$[1.8,1.9)$ & 0.40 & 0.60 & 0 & 526\\
$[2.0,2.1)$ & 0.84 & 0.16 & 0 & 476\\
$[2.1,2.2)$ & 0.78 & 0.10 & 0.12 & 455\\
$[2.2,2.3)$ & 0.26 & 0.10 & 0.64 & 432\\
$[2.4,2.5)$ & 0.04 & 0.04 & 0.92 & 395\\
$[2.8,2.9)$ & 0 & 0 & 1.00 & 336\\
\bottomrule
\end{tabular}
\caption{Types of attractor along $\gamma\tau$ for $\rho=25$, $\nu=0.005$, $p=0.05$.}
\label{Stab:tau}
\end{table}

For $\gamma\tau<2$ the attractor is quasi-periodic ($\lambda_1\approx0$) or a mode-locked periodic orbit ($\lambda_1<0$) whose period is close to the number of steps per cycle (e.g.\ 948, 920, 900, \dots, 661). The invariant closed curve that carries this motion exists as a circle homeomorphism only for $\gamma\tau\le1.39$ and for $1.54\le\gamma\tau\le1.74$ (Section~\ref{Ssec:circle}); in between, and up to $\gamma\tau\approx1.98$, the curve is folded although the attractor is not yet chaotic. The fraction of mode-locked values grows with $\gamma\tau$, up to $40$--$60\%$ for $\gamma\tau\approx1.4$--$1.9$. With the per-step threshold $\lambda_1>0.003$ used in Table~\ref{Stab:tau}, the first chaotic value appears at $\gamma\tau=2.18$, after which $\lambda_1$ grows continuously from zero ($0.003$ near $\gamma\tau=2.2$, $0.006$ at $2.5$, $0.01$ at $3.2$). Because a cycle lasts $450$--$950$ steps along this line, this threshold corresponds to $1.4$--$2.8$ per cycle and is conservative near the onset. Section~\ref{Ssec:circle} follows the route with exponents per cycle and identifies the breakdown of the invariant curve directly.

\subsection{Reduction to a circle map: torus breakdown and horseshoe}\label{Ssec:circle}

\emph{Definition.} For an orbit of $F_\tau$ let $n_k$ be the step at which $\log_{10}a$ first reaches $-3$ during the $k$-th outbreak, and $t_k=n_k-1+\theta_k$, $\theta_k\in[0,1)$, the crossing time obtained by linear interpolation of $\log_{10}a$ between steps $n_k-1$ and $n_k$. The outbreak phase is $\varphi_k=t_k\bmod 1$, and $\tilde g(\varphi_k)=\varphi_k+(t_{k+1}-t_k)$. Shifting an orbit by one step does not change $\varphi$, and on a one-dimensional attractor the interval $t_{k+1}-t_k$ is a function of $\varphi_k$; $\tilde g$ is then the lift of a degree-one circle map $g$, and its rotation number is the fractional part of the mean number of steps per outbreak. The phase map is shown in Fig.~\ref{fig:circlemap}(a).

\emph{Validity of the reduction.} The Lyapunov exponents $\lambda_1,\lambda_2,\lambda_3$ along this line are, per step, $0.000$, $-0.0106$, $-0.0375$ at $\gamma\tau=1.2$; $-0.0003$, $-0.0108$, $-0.0453$ at $1.6$; $0.0015$, $-0.0132$, $-0.0510$ at $2.0$; $0.0000$, $-0.0124$, $-0.0647$ at $2.5$ (periodic window); and $0.0080$, $-0.0214$, $-0.0707$ at $3.0$. Multiplied by the number of steps per outbreak, $\lambda_2$ gives a transverse contraction of $e^{-8.5}$, $e^{-6.6}$, $e^{-6.4}$, $e^{-4.8}$ and $e^{-6.8}$ per cycle, so that after one or two outbreaks any orbit lies on the attractor to high accuracy and the phase of the next outbreak is, to that accuracy, a function of the phase of the current one. We measured the residual scatter of the points $(\varphi_k,\tilde g(\varphi_k))$ about a local quadratic fit (25 neighbours). Its median and $99\%$ quantile are $3\times10^{-6}$ and $7\times10^{-4}$ at $\gamma\tau=1.8$, $4\times10^{-4}$ and $2\times10^{-3}$ at $2.0$, $8\times10^{-3}$ and $3\times10^{-2}$ at $2.5$, and $4\times10^{-2}$ and $0.11$ at $3.0$. The scatter grows with $\gamma\tau$ because the Kaplan--Yorke dimension of the attractor grows; at the reference point ($\rho=46$, $\gamma\tau=5.5$, $D_{KY}=1.81$) the reduction no longer holds, and the horseshoe there is studied through the unstable manifold of a saddle orbit (Section~\ref{Ssec:horseshoe}).

\emph{Loss of invertibility.} The derivative $\tilde g'(\varphi_k)$ was computed along orbits by complex-step differentiation of $F_\tau$ (step $10^{-20}$), which gives the derivative of the crossing time to machine precision, renormalising the tangent vector at each crossing. Sampling the curve with the transients of 300 initial conditions (three cycles after two cycles of transient) and taking the $0.5\%$ quantile as the minimum slope gives $0.90$ at $\gamma\tau=1.0$, $0.41$ at $1.2$, $0.14$ at $1.3$, $0.020$ at $1.38$ and $0.002$ at $1.39$; the slope becomes negative at $1.40$ ($-0.009$) and remains between $-0.01$ and $-0.05$ up to $1.51$, with larger negative values at $1.43$, $1.52$ and $1.53$. The map is invertible again for $1.54\le\gamma\tau\le1.74$ (minimum slope between $0.014$ and $0.29$) and non-invertible from $1.75$ on, with minimum slope $-0.49$ at $1.80$, $-3.2$ at $2.0$, $-16$ at $2.5$ and $-31$ at $3.0$ (Fig.~\ref{fig:circlemap}(b)). A circle homeomorphism has a unique rotation number, its orbits are quasi-periodic or converge to periodic orbits, and its topological entropy is zero~(de Melo \& van Strien, 1993, cited in the main text); hence there is no chaos for $\gamma\tau\le1.39$ or for $1.54\le\gamma\tau\le1.74$. A non-invertible degree-one circle map can have an interval of rotation numbers and positive entropy, and this is what happens here.

\emph{Exponents per cycle.} With the same complex-step derivative along a single orbit of $1500$ cycles, the Lyapunov exponent per cycle $\lambda_c=\langle\ln|\tilde g'|\rangle$ is zero (quasi-periodic) or negative (mode-locked) wherever $g$ is invertible. The first positive values are isolated: $0.17$ at $\gamma\tau=1.80$ and $0.46$ at $1.88$. From $1.98$ on $\lambda_c$ is positive except in periodic windows ($0.63$ at $1.98$, $1.0$ at $2.08$, $1.5$ at $2.30$, $2.1$ at $2.60$, $2.6$ at $3.0$; Fig.~\ref{fig:circlemap}(c)), and on the chaotic attractor about $40$--$50\%$ of the cycles have $\tilde g'<0$, i.e.\ pass through the fold.

\emph{Horseshoe.} Let $\mu=\max(0.05,\,5q_{99})$, where $q_{99}$ is the $99\%$ quantile of the scatter. On each monotone lap of $\tilde g$ we counted the integers $n$ such that the lap covers $[n-\mu,\,n+1+\mu]$; each such $n$ gives an arc of the circle whose image covers the whole circle with margin $\mu$, and arcs on the same or on different laps are disjoint. With $k$ such arcs, $g$ has a full $k$-horseshoe and its topological entropy is at least $\ln k$ per cycle. We find $k=3$ at isolated values from $\gamma\tau=2.16$ and for most values from $2.31$, $k=5$ from $2.45$, $k=7$ from $2.57$, $k=9$ from $2.68$ and $k=11$ from $2.82$ ($k=13$ at $2.99$). The bound $\ln k$ is close to the exponent per cycle: $\ln5=1.61$ against $\lambda_c=2.1$ at $2.52$, $\ln9=2.20$ against $2.21$ at $2.70$, and $\ln11=2.40$ against $2.61$ at $3.0$. Divided by the steps per outbreak, $\ln11$ gives at least $0.0076$ per step at $\gamma\tau=3$, against $\lambda_1=0.0080$.

\emph{Towards a computer-assisted proof.} Two ingredients of the argument above are numerical: the reduction to $g$ and the covering of the arcs. Both can in principle be validated with interval arithmetic by covering relations for h-sets~(Zgliczy\'nski \& Gidea, 2004, cited in the main text): each arc would be thickened into a thin box with one expanding coordinate (the phase) and two contracting ones, and one would check with rigorous enclosures of $F_\tau^n$ that the image of each box crosses every box in the expanding direction while staying inside them in the contracting directions. The conditions are favourable. The margins ($\mu=0.05$--$0.57$) exceed the scatter by a factor of five and the floating-point error by more than eight orders of magnitude: recomputing three successive crossing times in 80-bit arithmetic changed them by at most $7\times10^{-12}$ at $\gamma\tau=2.0$ and $2.3\times10^{-10}$ at $\gamma\tau=3.0$. The transverse contraction per cycle ($e^{-5}$--$e^{-9}$) makes the stable directions easy to control. The main technical difficulty is the length of the cycle (about $320$ steps at $\gamma\tau=3$), which requires an enclosure method that controls the wrapping effect, such as Lohner's algorithm; this is standard in existing libraries for rigorous computation.

%======================================================================
\section{Saddle periodic orbits}\label{Ssec:upo}

At the reference point ($\rho=46$, $\nu=0.002$, $p=0.05$, $\gamma\tau=5.5$) we searched for periodic orbits by Newton's method on $F_\tau^k(x)=x$, with finite-difference Jacobians, starting from near-recurrences of a long orbit on the attractor, for $k=40,\dots,129$. We found 31 orbits (residual $<10^{-14}$), all with one unstable and two stable multipliers:
\begin{itemize}
\item one-outbreak orbits for every period $k=58$--$79$, each consisting of one outbreak followed by a recharge; these are exactly the inter-outbreak intervals observed on the attractor (Section~\ref{Ssec:symb});
\item two-outbreak orbits for $k=120$--$129$ (except 124).
\end{itemize}
The unstable multipliers are $|\Lambda_u|\approx10$--$28$ for one-outbreak orbits and $\approx165$--$507$ for two-outbreak orbits, of the order of products of two one-outbreak multipliers, as expected for horseshoe-like dynamics. The sign of $\Lambda_u$ alternates between orbits (saddles with and without orientation reversal), as for the two branches of a horseshoe. The growth rates per step, $\ln|\Lambda_u|/k=0.03$--$0.05$, are comparable with $\lambda_1=0.037$ (Table~\ref{Stab:upo}).

\begin{table}[h]
\centering\small
\begin{tabular}{@{}rrrr@{\qquad}rrrr@{}}
\toprule
$k$ & $\Lambda_u$ & $\Lambda_{s,1}$ & $\ln|\Lambda_u|/k$ & $k$ & $\Lambda_u$ & $\Lambda_{s,1}$ & $\ln|\Lambda_u|/k$\\
\midrule
58 & 9.62 & 0.098 & 0.039 & 70 & 25.53 & 0.013 & 0.046\\
59 & 13.86 & 0.062 & 0.045 & 72 & 24.78 & 0.012 & 0.045\\
60 & 16.76 & 0.047 & 0.047 & 75 & $-25.02$ & $-0.028$ & 0.043\\
61 & $-17.27$ & $-0.065$ & 0.047 & 77 & 18.82 & 0.013 & 0.038\\
62 & 20.79 & 0.032 & 0.049 & 79 & $-11.12$ & $-0.044$ & 0.030\\
64 & $-23.77$ & $-0.045$ & 0.050 & 120 & $-233.5$ & $-0.0042$ & 0.045\\
66 & 24.89 & 0.019 & 0.049 & 122 & $-323.3$ & $-0.0028$ & 0.047\\
68 & $-28.42$ & $-0.034$ & 0.049 & 127 & $-506.5$ & $-0.0012$ & 0.049\\
\bottomrule
\end{tabular}
\caption{Selected saddle periodic orbits at the reference point. $\Lambda_u$: unstable multiplier; $\Lambda_{s,1}$: largest stable multiplier (the third is of order $10^{-8}$ for $k\le79$). One-outbreak orbits exist for every $k=58$--$79$.}
\label{Stab:upo}
\end{table}

%======================================================================
\section{Unstable manifold and horseshoe}\label{Ssec:horseshoe}

\begin{figure}[h]
\centering
\includegraphics[width=\linewidth]{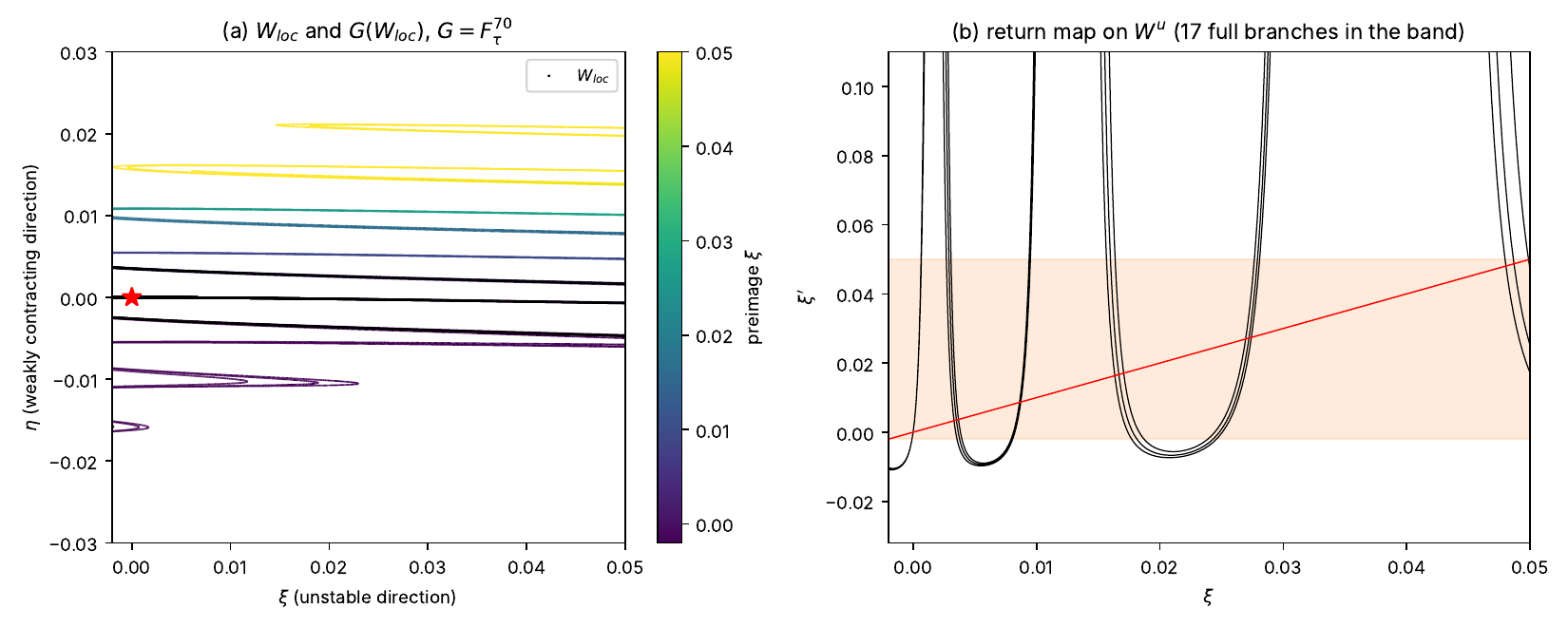}
\caption{(a) Local unstable manifold $W_{loc}$ of the period-70 saddle (black; star: point of the orbit) and its image under $G=F_\tau^{70}$, in local coordinates $(\xi,\eta)$, coloured by the preimage $\xi$: the image returns to the section in several folded layers. (b) Return map $\xi\mapsto\xi'$ on $W^u$; the shaded band contains 17 full branches.}
\label{Sfig:horseshoe}
\end{figure}

We took the period-70 saddle ($\Lambda_u=25.53$), with a point in the recharge phase, and iterated a fundamental segment of its local unstable manifold (initial distances from $10^{-8}$, four iterations of $F_\tau^{70}$). Figure~\ref{Sfig:horseshoe} shows the manifold on a section of the recharge phase in local coordinates $(\xi,\eta)$ (unstable and weakly contracting directions) and the induced one-dimensional return map $\xi\mapsto\xi'$.

The manifold is stretched by a factor of about $25$ per turn and folded: its image returns to the section in several nearly parallel layers. Stretching takes place along the recharge phase and folding during the outbreak. The return map has monotone branches that cover the whole interval: 17 in the band shown in Fig.~\ref{Sfig:horseshoe}(b), and between 7 and 20 depending on the interval and covering criterion used. If these branches define a topological horseshoe, the topological entropy of $F_\tau$ is bounded below by $\ln7/70=0.028$ to $\ln20/70=0.043$ per step. This bound is conditional, since the covering has been verified numerically and not with rigorous (interval) arithmetic.

%======================================================================
\section{Symbolic dynamics of inter-outbreak intervals}\label{Ssec:symb}

\begin{figure}[h]
\centering
\includegraphics[width=\linewidth]{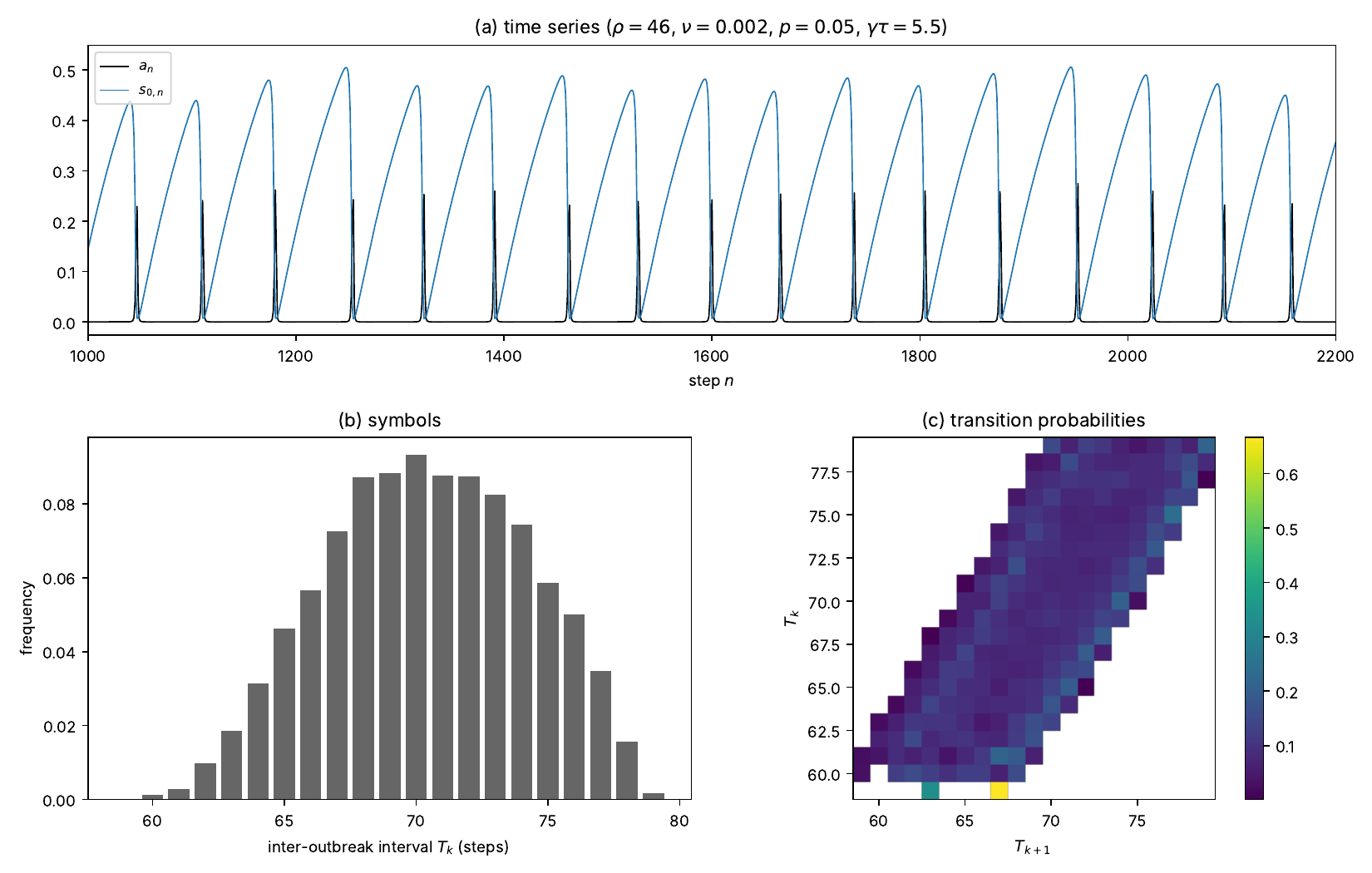}
\caption{(a) Time series of $a_n$ and $s_{0,n}$ at the reference point. (b) Frequencies of the 21 symbols (inter-outbreak intervals). (c) Transition probabilities $T_k\to T_{k+1}$; white cells are transitions not observed.}
\label{Sfig:symb}
\end{figure}

The inter-outbreak intervals of an orbit with 28\,403 outbreaks (upward crossings of $a=0.05$) take 21 integer values, from 59 to 79 steps, which we use as symbols (Fig.~\ref{Sfig:symb}).
\begin{itemize}
\item \emph{Grammar.} Of the $21^2=441$ possible ordered pairs of consecutive symbols, 216 (49\%) occur. Allowed transitions form a band, $T_{k+1}$ lying within a few steps of $T_k$, so that the interval drifts over several outbreaks rather than jumping between extremes. The only highly probable transitions start from the shortest interval, which is followed almost always by 63 or 67 steps.
\item \emph{Entropy.} The topological entropy of this first-order grammar (logarithm of the spectral radius of the adjacency matrix), divided by the mean interval, is about $0.031$ per step. It is consistent with the horseshoe bound ($0.028$--$0.043$) and below $\lambda_1=0.037$.
\item \emph{Interval map.} The relation $T_k\mapsto T_{k+1}$ is not a function: each symbol has up to 11 successors (2 for the shortest interval), because the state at the start of a trough depends on more than one variable. The dynamics does not reduce to a one-dimensional map of intervals.
\end{itemize}

%======================================================================
\section{Correlation dimension}\label{Ssec:d2}

\begin{figure}[h]
\centering
\includegraphics[width=0.85\linewidth]{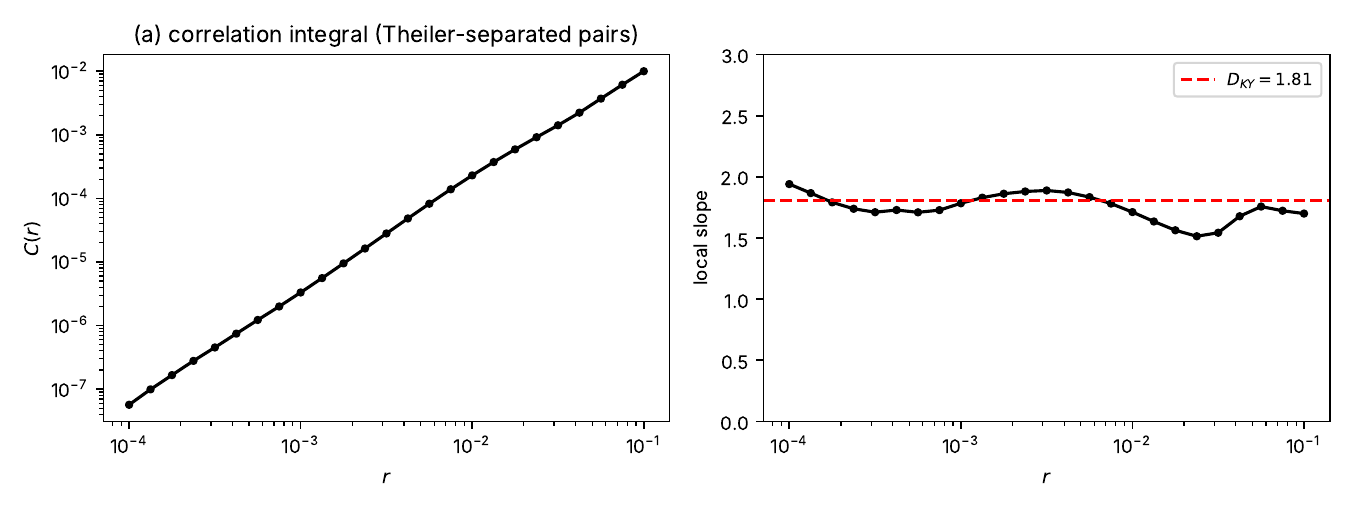}
\caption{Correlation integral with Theiler-separated pairs (left) and local slope (right) at the reference point. Dashed line: $D_{KY}=1.81$.}
\label{Sfig:d2}
\end{figure}

We computed the Grassberger--Procaccia correlation integral at the reference point, in coordinates $(s_0,s_1,\log_{10}a)$ normalised to zero mean and unit variance. We took two segments of $10^5$ points each (every fourth step of an orbit of $2\times10^6$ steps), separated by $4\times10^5$ steps, and counted only pairs with one point in each segment, so that all pairs are separated in time by at least $4\times10^5$ steps (Theiler window). $C(r)$ was evaluated at 25 logarithmically spaced radii in $[10^{-4},10^{-1}]$. The scaling range was chosen as the interval in which the local slope is approximately constant, $r\in[10^{-4},10^{-2}]$; above $r\approx10^{-2}$ the slope decreases to about 1.5 because of the finite size of the attractor. Least-squares fits over the whole scaling range give $D_2=1.80$, and fits over seven sub-windows of one to two decades within $[10^{-4},10^{-1.5}]$ give values between $1.75$ and $1.85$ (Table~\ref{Stab:d2fits}, Fig.~\ref{Sfig:d2}). This agrees with the Kaplan--Yorke dimension $D_{KY}=1.81$ obtained from the Lyapunov spectrum, and confirms that the attractor is a fractal set of dimension between one and two, a thin band folded infinitely many times.

\begin{table}[h]
\centering\footnotesize
\setlength{\tabcolsep}{4pt}
\begin{tabular}{@{}lccccccc@{}}
\toprule
window in $\log_{10}r$ & $[-4,-2]$ & $[-4,-3]$ & $[-3.5,-2.5]$ & $[-3,-2]$ & $[-3.5,-2]$ & $[-4,-2.5]$ & $[-3,-1.5]$\\
\midrule
$D_2$ & 1.80 & 1.75 & 1.79 & 1.85 & 1.82 & 1.77 & 1.78\\
\bottomrule
\end{tabular}
\caption{Correlation dimension at the reference point for different fitting windows.}
\label{Stab:d2fits}
\end{table}

\end{document}